\documentclass[final]{siamonline171218}
\usepackage{amsmath}
\usepackage{amsfonts}
\usepackage{amssymb}
\usepackage{color}
\usepackage{enumerate}

\usepackage{graphicx} % Allows including images
\usepackage{booktabs}
\usepackage[caption=false]{subfig}
\usepackage{cleveref}
\usepackage[utf8]{inputenc}
\usepackage[english]{babel}
\usepackage{algorithmic}
\Crefname{ALC@unique}{Line}{Lines}

\usepackage{lipsum}
\usepackage{amsfonts}
\usepackage{graphicx}
\usepackage{epstopdf}
\usepackage{algorithmic}
\ifpdf
  \DeclareGraphicsExtensions{.eps,.pdf,.png,.jpg}
\else
  \DeclareGraphicsExtensions{.eps}
\fi

\numberwithin{theorem}{section}

\newcommand{\TheTitle}{Continuation with Non-invasive Control Schemes: Revealing Unstable States in a Pedestrian Evacuation Scenario} 
\newcommand{\TheAuthors}{I. Panagiotopoulos, J. Starke, J. Sieber, W. Just}

\headers{Non-invasive Control of a Particle Model}{\TheAuthors}

\title{{\TheTitle}\thanks{Submitted to the editors DATE.
}}

\author{
  Ilias Panagiotopoulos\thanks{Institute of Mathematics, University of Rostock, Germany.}
  \and
  Jens Starke\thanks{Institute of Mathematics, University of Rostock, Germany.}
  \and
  Jan Sieber\thanks{College of Engineering, Mathematics and Physical Sciences, University of Exeter, United Kingdom.}
      \and
  Wolfram Just\thanks{School of Mathematical Science, Queen Mary University of London, United Kingdom.}
}

\usepackage{amsopn,mathtools}

\DeclareMathOperator{\tr}{tr}
\DeclareMathOperator{\spec}{spec}
\DeclareMathOperator{\re}{Re}

\DeclareMathOperator{\std}{std}
\DeclareMathOperator{\tol}{tol}

\newcommand{\tran}{\mathsf{T}}
\newcommand{\R}{\mathbb{R}}
\newtheorem{assumption}{Assumption}[section]
\newcommand{\js}[1]{{#1}}
\ifpdf
\hypersetup{
  pdftitle={\TheTitle},
  pdfauthor={\TheAuthors}
}
\fi

\begin{document}

\maketitle

% REQUIRED
\begin{abstract}

  This paper presents a framework to perform bifurcation
  analysis in laboratory experiments or simulations.  We employ
  \emph{control-based continuation} to study the dynamics of a
  macroscopic variable of a microscopically defined model, exploring %testing
  the
  potential viability of the underlying feedback control techniques in
  an experiment.  In contrast to previous experimental studies that used iterative
  root-finding methods on the feedback control targets, we propose a
  feedback control law that is inherently non-invasive. That is, the
  control discovers the location of equilibria and stabilizes them
  simultaneously. We call the proposed control \emph{zero-in-equilibrium feedback control} and we prove that
  it is able to stabilize branches
  of equilibria, except at singularities of codimension $n+1$, where $n$ is the number of state space dimensions the
  feedback can depend on.

  We apply the method to a simulated evacuation scenario were pedestrians have
  to reach an exit after maneuvering left or right around an obstacle.
  The scenario shows a hysteresis phenomenon with bistability and
  tipping between two possible steady pedestrian flows in microscopic
  simulations.  
  We demonstrate for the evacuation scenario that the proposed
  control law is able to uniformly discover and stabilize steady
  flows along the entire branch, including points where other
  non-invasive approaches to feedback control become singular.

\end{abstract}

% REQUIRED

% REQUIRED
\section{Introduction}

The analysis of systems involving many interacting microscopic
components is often desired in terms of a few, macroscopic
quantities, such as averages over all components
\cite{givon2004}. Qualitative behavior, such as the system being in
equilibrium or whether it shows multistability or is near a tipping
point, are expressed at this macroscopic level
\cite{Scheffer2009a}. 
However, the derivation of a model for the evolution of the macroscopic
quantities 
typically 
relies on assumptions that are not realistic or are known to introduce a
bias. For example in networks, mean-field equations  rely on closure
approximations. These closures assume absence of correlations beyond a
fixed diameter, for practical reasons this diameter equals $1$ such
that one ignores correlations beyond nearest neighbors
\cite{GDB06,GK08,Pellis2015,Kiss2017}.
On the other hand, a microscopic model may well be amenable
to direct simulations (from which one can extract macroscopic
quantities) and may be easy to connect to observed data or first
principles, as its parameters encode the individual behavior of the
interacting agents of the underlying system. Applications of direct
simulations for individual-based models in ecology are discussed in
\cite{railsback2019agent}.

We focus on an approach inspired by its applicability to experiments\js{, providing a non-technical overview of in \Cref{sec:bg+results}}.
%To overcome this problem, 
The engineering and physics communities have independently developed
feedback control laws that enable one to perform bifurcation analysis
directly on physical experiments \cite{PYRAGAS1992421,WANG19951213,barton2017control}. In particular, these control-based methods can track dynamical
phenomena that are either dynamically unstable or too sensitive to
disturbances to be visible in uncontrolled experiments \cite{PhysRevLett.100.244101}. \Cref{sec:control-based} will start with a brief review of methods for performing bifurcation analysis without deriving explicit equations, including our approach of using \emph{non-invasive feedback control}. 

\Cref{sec:inh:noninv} will reformulate one of the classical inherently non-invasive feedback laws, which are based on washout filters, in a way that makes it treatable with standard controllability arguments, which can be decided by the regularity or singularity of a controllability matrix. We then construct a new control law that removes all singularities, except for some events of high codimension along branches of equilibria.

We demonstrate the proposed methodology on a microscopically defined model, a particle flow model for pedestrians moving along a corridor past an obstacle, introduced in \Cref{sec:scenario}. Our control-based continuation of this scenario, which exhibits bistability  and two tipping points at a macroscopic level, reveals the unstable pedestrian flows and completes a bifurcation diagram in \Cref{sec:cbped} without the use of any macroscopic model. During our analysis, we assume all the limitations of a physical experiment and, thus, our approach can be extended to real life scenarios.

\section{Non-technical overview}
\label{sec:bg+results}

Let us assume that an experiment (computational or physical) can be described by an ordinary differential equation (ODE) with some state $x(t)$ and parameters $\mu$ from which the output originates in the form $y(t)=g(x(t))$. \emph{Feedback control} takes the output $y(t)$ and feeds back a control input signal $u(t)$, which depends on $y$ and possibly its history. In this section we discuss the case of state feedback
control, $y(t)=x(t)$, to simplify notation.
Feedback control needs to be designed in a way to be stabilizing and, in addition for the purposes discussed in our paper, non-invasive.
 Making feedback control
stabilizing is a standard control theoretical task involving the
construction of state observers (if necessary) and \emph{control
  gains}, which are amplification factors for $y$ entering the input
$u$. Non-invasiveness refers to the property that the control input
signal $u(t)$ vanishes in the stabilized steady states after transients have settled. When $u=0$, up
to disturbances typical for experiments, then one observes phenomena
of the original uncontrolled system, where $u$ was $0$.

\subsection{Design of inherently non-invasive feedback control}
\label{sec:intro:design}
The most well-known example of non-invasive feedback control is
time-delayed feedback \cite{PYRAGAS1992421}, where
$u(t)=K\cdot(x(t)-x(t-T))$.  This input automatically
vanishes whenever the $x(t)$ settles to an equilibrium or a periodic
orbit of period $T$, but is not able to stabilize equilibria or
periodic orbits of forced systems with single eigenvalues $0$ or $1$, respectively
\cite{hovel2010control}. For continuation of equilibria two other
classes of non-invasive control laws have been designed and
investigated. First, washout filters
\cite{AWC94,hassouneh2004washout} % assume
% that the state is governed by an ordinary differential equation (ODE)
% with input $u$ and
add extra degrees of freedom, $x_\mathrm{wo}$, which we
%formulate feedback control with washout filter 
 formulate in the form ($x\in\R^{n_x}$ and $u\in\R^{n_u}$)
\begin{align}
  \label{intro:washout}
  \dot x(t)&=f(x(t),u(t))\mbox{,}& \dot x_\mathrm{wo}(t)&=u(t)\mbox{,}
  &u(t)=K_\mathrm{st}[x(t)-x_\mathrm{ref}]+K_\mathrm{wo}[x_\mathrm{wo}(t)-x_\mathrm{wo,ref}]\mbox{,}
\end{align}
where $x_\mathrm{ref},x_\mathrm{wo,ref}$ are constant reference values chosen by the experimenter.
The feedback is non-invasive as any equilibrium  of
\eqref{intro:washout} is also an equilibrium $x_\mathrm{eq}$ of
$\dot x=f(x,0)$. One can find gains $(K_\mathrm{st},K_\mathrm{wo})$ to
make $x_\mathrm{eq}$ stable with arbitrary decay rate if and only if
the matrix $\partial_1f(x_\mathrm{eq},0)$ is regular and the pair
$(\partial_1f(x_\mathrm{eq},0),\partial_2f(x_\mathrm{eq},0))$ is controllable. Formulation \eqref{intro:washout} is \js{a generalization of the construction in} the original papers \cite{AWC94,hassouneh2004washout}, showing
that the additional degrees of freedom form the integral component
of % , when considering \eqref{intro:SP} as
a proportional-integral (PI) control, which enforces $u=0$ in the
equilibrium. The PI control in \eqref{intro:washout} is degenerate, as the right-hand side for
$\dot x_\mathrm{wo}$ depends on $x$ only through $u$. This simplifies the
proof of the stabilization criterion compared to the original sources
\cite{bazanella1997control,hassouneh2004washout} and reduces the
design of the feedback gains $(K_\mathrm{st},K_\mathrm{wo})$ to a
standard linear control design problem (see \Cref{thm:washout}). % \js{In \cref{sec:SO} we will also show that the most general formulation of washout filters in \cite{hassouneh2004washout} is a special case of framework \eqref{intro:washout}.}

When continuing branches of equilibria, control law \eqref{intro:washout}
fails at saddle-node (fold) bifurcations, where
$\partial_1f(x_\mathrm{eq},\mu_\mathrm{eq})$ is singular, motivating a second type of non-invasive feedback law. 
Assuming that the dynamical system depends
on a parameter $\mu\in\R$ and has a branch of equilibria,
parameterized by a scalar,
$s\mapsto(x_\mathrm{eq}(s),\mu_\mathrm{eq}(s))$, 
Siettos
\emph{et al}.~\cite{siettos2004coarse} proposed to dynamically adjust
the parameter $\mu$ using the control input $u$:
\begin{align}
  \label{intro:SP}
  \dot x(t)&=f(x(t),\mu(t))\mbox{,}& \dot \mu(t)&=u(t)\mbox{,}
  &u(t)=K_{\mathrm{st},x}[x(t)-x_\mathrm{ref}]+K_{\mathrm{st},\mu}[\mu(t)-\mu_\mathrm{ref}]\mbox{.}
\end{align}
This feedback control law is also non-invasive in the above sense,
such that one may track the branch of equilibria using
pseudo-arclength continuation
\cite{allgower2012numerical,DS13,doedel1997continuation,govaerts2000numerical,K04}
with a sequence of points $(x_\mathrm{ref},\mu_\mathrm{ref})$
predicted by the pseudo-arclength continuation algorithm. 
A limiting case of \eqref{intro:SP} was used in experiments in
\cite{BartonSieber2013}. One can find gains
$(K_{\mathrm{st},x},K_{\mathrm{st},\mu})$ to make
$(x_\mathrm{eq},\mu_\mathrm{eq})$ stable with arbitrary decay rate if
and only if the pair
$(\partial_1f(x_\mathrm{eq},\mu_\mathrm{eq}),\partial_2f(x_\mathrm{eq},\mu_\mathrm{eq}))$
is controllable (see \Cref{thm:SP}). In contrast to
\eqref{intro:washout}, control law \eqref{intro:SP} does not fail near saddle-node
bifurcations, which are points of particular interest when performing
bifurcation analysis. However, \Cref{thm:SP} implies that controllability through a scalar
bifurcation parameter $\mu$ will break down at codimension-$1$ events
such that one may generically encounter singularities in
single-parameter continuations.
\paragraph{Zero-in-equilibrium feedback control}
In this paper we generalize \eqref{intro:washout} and \eqref{intro:SP}
to introduce a control law that avoids these singularities by
combining \eqref{intro:washout} and \eqref{intro:SP}. We formulate the
law for the case where the ODE depends on a scalar bifurcation
parameter $\mu$ and an additional scalar control input $u\in\R$. This
permits us to introduce an additional scalar gain $a\in\R$ to propose
the feedback
\begin{align}
  \label{intro:combi}
  \dot x(t)&=f(x(t),\mu(t),au(t))\mbox{,}& \dot \mu(t)&=u(t)\mbox{,}
  &u(t)=K_{\mathrm{st},x}[x(t)-x_\mathrm{ref}]+K_{\mathrm{st},\mu}[\mu(t)-\mu_\mathrm{ref}]\mbox{.}
\end{align}
For \eqref{intro:combi} one can find gains
$(a,K_{\mathrm{st},x},K_{\mathrm{st},\mu})$ to make
$(x_\mathrm{eq},\mu_\mathrm{eq})$ stable with arbitrary decay rate if
and only if the matrix pair $(R_\mu,f_xR_u)$ is regular (polynomial
$\lambda\mapsto\det(R_\mu+\lambda f_xR_u)\not\equiv0$, see
\Cref{thm:combi}), where
$f_x=\partial_1f(x_\mathrm{eq},\mu_\mathrm{eq},0)$, and $R_\mu$ and
$R_u$ are the controllability matrices of $f_x$ with respect to
$f_\mu=\partial_2f(x_\mathrm{eq},\mu_\mathrm{eq},0)$ and
$f_u=\partial_3f(x_\mathrm{eq},\mu_\mathrm{eq},0)$, respectively: for  $n_x=\dim x$
\begin{align}
  \label{intro:combi:cond}
  R_u&=[f_u,f_xf_u,\ldots,f_x^{n_x-1}f_u]\in\R^{n_x\times n_x}\mbox{,}&
  R_\mu&=[f_\mu,f_xf_\mu,\ldots,f_x^{n_x-1}f_\mu]\in\R^{n_x\times n_x}\mbox{.}  
\end{align}
This regularity condition is violated only in events of codimension
$n_x+1$.  Thus, \eqref{intro:combi} with suitable gains stabilizes the
natural equilibria of an ODE dynamically uniformly along the whole
branch for generic branches of equilibria.  This is in contrast to
\eqref{intro:washout} and \eqref{intro:SP}, which one expects to fail
at isolated points of the curve (codimension-$1$ events), namely when
$f_x$ is singular (for \eqref{intro:washout}, assuming that $R_u$ is
always regular, as we can freely choose a suitable $u$), or when
$R_\mu$ is singular (for \eqref{intro:SP}). 

We use the term \emph{zero-in-equilibrium feedback control} to indicate that the class of feedback laws \eqref{intro:combi} contains washout filters, \eqref{intro:washout}, and control through parameter, \eqref{intro:SP}, as limiting cases ($a\to\infty$ and $a=0$), but is less general than ``non-invasive'' which includes time-delayed feedback. 

\js{\Cref{fig:flow_manipulation_js} in \cref{sec:control-based} shows a sketch how control through the bifurcation parameter, \eqref{intro:SP} and zero-in-equilibrium control \eqref{intro:combi} affect the flow near the equilibrium branch in the case of scalar $x$.} 
\subsection{Demonstration on multi-particle model for pedestrian flow}
\label{sec:intro:pedestrian}
\begin{figure} [ht]
  \centering
  \includegraphics[width=\textwidth]{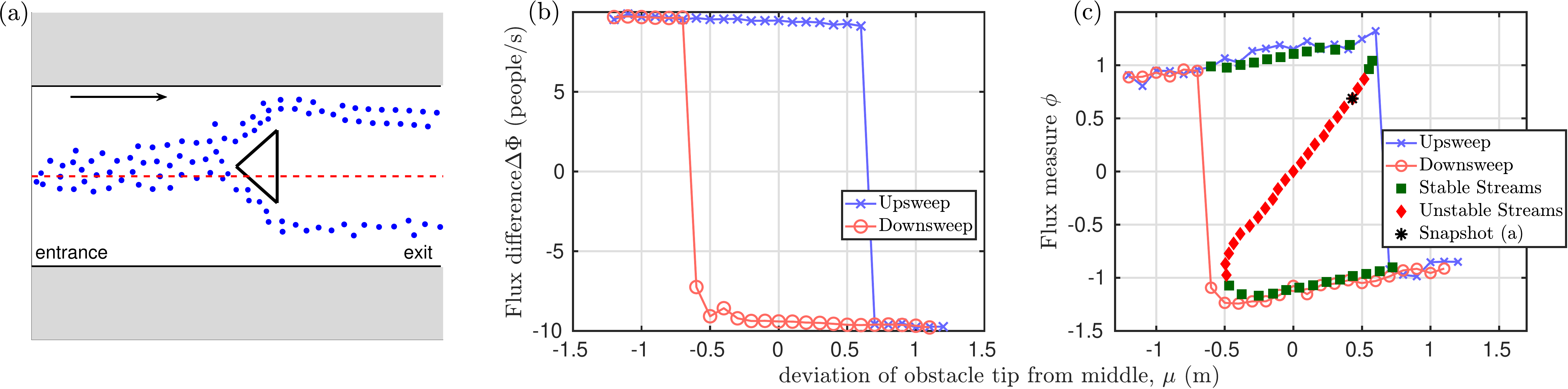}
  \caption{Response of a flow of $N=100$ pedestrians through a
    corridor to an obstacle, depending on the obstacle position $\mu$,
    relative to middle of corridor. \textup{(a)} Top view snapshot of
    pedestrians \textup{(}blue dots\textup{)} moving through a
    $20\,\mathrm{m}\times10\,\mathrm{m}$ corridor past an
    obstacle. \textup{(}b\textup{)} Parameter sweep, changing distance $\mu$ of 
    obstacle tip from red dashed line in \Cref{fig:intro}\textup{(a)}, and response $\Delta\Phi$, given in \eqref{eq:DPhi}. The
    coordinates of the snapshot \Cref{fig:intro}\textup{(a)} are
    indicated as a star symbol in \Cref{fig:intro}\textup{(c)}. \textup{(c)} Full bifurcation diagram
     obtained using the newly proposed control law \eqref{intro:combi} and the spatially averaged flux measure $\phi$, given in
    \eqref{flux}. % To facilitate the analysis, an spatially averaged flux of pedestrians $\phi$ is defined, see also \eqref{flux}.
  }
\label{fig:intro}
\end{figure}
To illustrate continuation with non-invasive feedback control, we
consider a particle model describing pedestrians moving through a
corridor past an obstacle in an evacuation scenario.

The microscopic behavior of every pedestrian (treated as a point
particle) is based on the social force model by Helbing and Molnar~\cite{HelbingSocial} with an additional preference of
pedestrians toward alignment with near-by others moving in roughly the
same direction \cite{starkeieeeqingdao14}. All pedestrians are moving towards the end of the corridor, as shown in
\Cref{fig:intro}(a).
Inside the corridor there is a triangular obstacle which blocks the
pedestrians' straight path to the exit. As a result, they have to
choose a route, left or right of the obstacle from their point of view
(see top view \cref{fig:intro}(a)).  The position $\mu$ of the
obstacle relative to the middle of the corridor (vertical distance of
triangle tip from red dashed line in \Cref{fig:intro}(a)) changes the
pedestrians' preference for each route. We consider this position
(measured in meters) as the system parameter $\mu$.  The macroscopic
variable of interest is the difference between the flows of
pedestrians along the two different routes. \Cref{fig:intro}(b) shows
the time-averaged flux difference $\Delta\Phi$, given in \eqref{eq:DPhi}, between left and right
route, measured in people per second.

This particle system exhibits a bistability and hysteresis phenomenon:
once the majority of pedestrians has chosen a particular
route, even a small alignment effect will cause pedestrians to follow
the flow, even if the current route is less direct than the
alternative, until there is a sudden transition of the pedestrian flow
to the other route. \Cref{fig:intro}(b) shows this effect in a
parameter study for obstacle position $\mu$, where the obstacle is
gradually shifted, first from $\mu=-1.25$\,m upwards to $1.25$\,m (in
blue, with crosses), then downwards again (in red, with
circles). There is a large region of bistability, which leads to the
hypothesis that the sudden transitions at macroscopic level are
saddle-node bifurcations, and that there is a branch of unstable
steady flows inside this bistability region, which acts as a threshold
for disturbances to cause spontaneous transition.

Control-based continuation using control law \eqref{intro:combi} confirms this hypothesis for the particle model as shown in \Cref{fig:intro}(c).
%discover this unstable branch and the saddle-node bifurcations by
Here the $y$-axis is a spatially averaged flux difference $\phi$, which the feedback control input $u$ depends on (see \eqref{flux} and \Cref{fig:pedweight} for definition of $\phi$). 
 The control-based continuation enables us to produce the full
bifurcation diagram, including both stable and unstable steady flows,
without the need of having access to any effective macroscopic
model. Control law \eqref{intro:combi} turns out to be especially
advantageous compared to parameter control \eqref{intro:SP} (which
would in principle also be feasible) because we are free to choose the control
input $u$ without additional computational cost. In laboratory experiments additional inputs may require
additional actuation equipment. We choose a bias force acting on
pedestrians directly in front of the obstacle (see \eqref{model:xddot+input} and \Cref{fig:box} for definition of input). The modulus
of $u$ is always small when controlling perturbations (due to random
entry of the pedestrians into the corridor) with this bias force, while \eqref{intro:SP}
caused large corrections in $\mu(t)$.

\section{Equation-free bifurcation analysis and control-based continuation}
\label{sec:control-based}
Various methods have been proposed to avoid the need for an explicit
macroscopic model for the bifurcation analysis required for high-level
qualitative analysis. Equation-free methods pioneered by Kevrekidis
\emph{et al}.~\cite{eqfreeKev,Kevrekidis:2010} %of them
have been
applied primarily to computational experiments originating from
multi-particle simulations, while methods based on feedback control
were developed for physical experiments
\cite{PYRAGAS1992421,WANG19951213,barton2017control}. 

This section reviews the two fundamentally different methods briefly. We
then show how unifying the known inherently non-invasive feedback
control for equilibria of nonlinear systems laws as special cases of
control with an integral component allows us to 
design control law \eqref{intro:combi}, which does not suffer from
singularities one would encounter along a generic branch of
equilibria. This makes \eqref{intro:combi} applicable to the particle
flow simulation of a pedestrian evacuation scenario such that we can
perform control-based continuation using inherently non-invasive
feedback control, without requiring numerical root-finding algorithms. Continuation using \eqref{intro:combi} also does not rely on information about partial derivatives at every step.

One conclusion
from our paper is that feedback control-based methods may also be an
easy-to-implement approach to equation-free analysis in computational
experiments.

\subsection{Equation-free analysis based on lift-evolve-restrict cycles}
\label{sec:eqfree}
The methodology originally proposed by Kevrekidis \emph{et al.}  (see e.g.,
\cite{eqfreeKev,Kevrekidis:2010} for reviews), named
\emph{equa\-tion-free analysis}, is able to perform high-level tasks,
such as bifurcation analysis or optimization of macroscopic behavior,
on complex simulations (such as multi-particle models) by judiciously
initialized simulations, without explicitly deriving a macroscopic
model. This is done by extracting numerical information about the
macroscopic behavior using suitable short simulation bursts of the
microscopic model as part of a \emph{lift-evolve-restrict} loop.  The
basic method requires the user to specify a projection of the
microscopic state into the space of macroscopic quantities of suitable
dimension (the restriction operator), and an embedding of the
(low-dimensional) macroscopic state space into the state space of the
microscopic complex simulation (the lifting operator). With the use of
data analysis techniques such as diffusion maps it may be possible to
determine the dimension and numerically optimal restriction projection
for the macroscopic variables automatically
\cite{chin2021enabling,coifman2008diffusion,singer2009detecting}. 
Equation-free analysis assumes that the \emph{lifting} operator maps
close to an assumed-to-exist attracting slow manifold, on which the
macroscopic evolution takes place.  To compensate for the error of not
lifting exactly on this manifold, an implicit formulation of the
lifting operator can be used
\cite{marschler2014implicit,SMS18,Vandekerckhove2011}.

\js{Equation-free analysis based on lift-evolve-restrict cycles has also been applied to agent-based models such as epidemic networks by Gross and Kevrekidis \cite{GK08} and a population of traders participating in a financial market by Siettos \emph{et al}.\ \cite{siettos2012equation} and Tsoumanis and Siettos \cite{tsoumanis2012detection}. The analysis in \cite{siettos2012equation} concluded by applying a washout filter similar to \eqref{intro:washout} to make an unstable steady state visible in simulations without lifting, while \cite{tsoumanis2012detection} applied systematic corrections to the bifurcation parameter similar to \eqref{intro:SP}.}
\subsection{Control-based continuation}
\label{sec:gendesign}
Mechanical engineering and physics research into discovering
dynamically unstable phenomena in physical experiments with nonlinear
behavior took a different approach to equation-free analysis, which is
more suitable to physical experiments
\cite{barton2017control,PYRAGAS1992421,PhysRevLett.100.244101,WANG19951213} as most experiments cannot be initialized at arbitrary points in state space. 

\subsubsection{Existence of underlying ODE and equilibria}
\label{sec:ode}

Control-based continuation assumes that the underlying dynamical system (which will be a stochastic multi-particle model
in our case) is governed by a
system of ordinary differential equations with a state $x(t)$
depending on a scalar parameter $\mu$, and with control inputs $u(t)$
and outputs $y(t)=g(x(t))$. This approach also makes assumptions concerning existence
of equilibria, and controllability and observability near these
equilibria as listed below.

In order to perform control-based continuation to an experiment, the control is applied with the feedback laws designed based on the aforementioned 
assumptions. Then one validates during the experiments whether the
feedback controlled system converges to an equilibrium to a
sufficiently good approximation given by the tolerances of the
experimental measurement equipment and expected disturbances.

Equivalently, when applying these experimental techniques to a
dynamical system given in the form of a simulation
(such as our particle model for pedestrians), then the
simulation is treated like a computational experiment. For systems with a large
number $N$ of interacting particles one may repeat the computational
experiment with different $N$ to observe whether a law of large
numbers holds, such that one has convergence of the feedback
controlled system to an equilibrium for increasing $N$.
\begin{assumption}[Equilibrium branch for system of ODEs]\label{ass:io}
  The dynamical system is assumed to be governed by a system of ODEs
  of the form
\begin{align}
  \label{eq:io:ode} \dot x(t)&=f(x(t),\mu,u(t))\mbox{,}& \mbox{where\ }f&:\R^{n_x}\times\R\times\R^{n_u}\to\R^{n_x}\mbox{.}
\end{align}
We assume that for $u(t)=0$, \eqref{eq:io:ode} has an isolated branch \textup{(}curve\textup{)} of
equilibria, parameterized by
$s\in[s_\mathrm{min},s_\mathrm{max}]\subset\R$,
$(x_\mathrm{eq}(s),\mu_\mathrm{eq}(s))$, %with the corresponding $y_\mathrm{eq}(s)=g(x_\mathrm{eq}(s),\mu_\mathrm{eq}(s))$,
satisfying
$0=f(x_\mathrm{eq}(s),\mu_\mathrm{eq}(s),0)=0$ for all
$s\in[s_\mathrm{min},s_\mathrm{max}]$.  
\end{assumption}
Initially we will discuss state feedback control, where we assume that
the control input $u$ may depend on all components of the state
$x(t)\in\R^{n_x}$. The for physical experiments more realistic case of
output feedback control, where only some output $y=g(x(t))$ may enter
the algebraic or dynamic rule for $u$ will be briefly discussed
afterwards.  The general criteria for non-invasiveness are identical
for output feedback to those of state feedback control. We will
consider output feedback in detail for the concrete criteria for
choosing feedback gains in \Cref{sec:1dbranch}. In a computational
experiment the full state is available, and one typically chooses a
few problem-specific quantities that enter the control
input. \Cref{fig:intro} showed two different flux measures,
$\Delta\Phi$ and $\phi$, as possible outputs.
\subsubsection{Controllability}
\label{sec:classical}
We abbreviate the partial derivatives of $f$ in the equilibria by
\begin{align*}
  f_x&:=\partial_1f(x_\mathrm{eq},\mu_\mathrm{eq},0)\in\R^{n_x\times n_x}\mbox{,}&
  f_\mu&:=\partial_2f(x_\mathrm{eq},\mu_\mathrm{eq},0)\in\R^{n_x\times 1}\mbox{,}\\
%  g_x&:=\partial_1 g(x_\mathrm{eq},\mu_\mathrm{eq})\phantom{,0}\in\R^{n_y\times n_x}\mbox{,}&
  f_u&:=\partial_3 f(x_\mathrm{eq},\mu_\mathrm{eq},0)\in\R^{n_x\times n_u}\mbox{}
%  g_\mu&:=\partial_1 g(x_\mathrm{eq},\mu_\mathrm{eq})\phantom{,0}\in\R^{n_y\times1}\mbox{}
\end{align*}
(dropping the argument $s$ in all expressions here). We recall that a
linear constant-coefficient system $\dot x=Ax+Bu$ %, $y=Cx$
with matrices $A\in\R^{n_x\times n_x}$ and $B\in\R^{n_x\times
  n_u}$ %and $C\in\R^{n_y\times n_x}$
is controllable if the matrix $[B,AB,\ldots,A^{n_x-1}B]$ has
full rank $n_x$ (for standard textbooks on control theory and controllability see \cite{aastrom2010feedback,richard2008modern}).
% , and it is linearly observable if
% $[C^\mathrm{T},A^\mathrm{T}C^\mathrm{T},\ldots,(A^\mathrm{T})^{n_x-1}C^\mathrm{T}]$
% has full rank $n_x$.
Controllability implies that there exist \emph{feedback gains}
$K_\mathrm{cn}\in\R^{n_u\times n_x}$ such that $A+BK_\mathrm{cn}$ is a
Hurwitz matrix. More precisely, the spectrum of $A+BK_\mathrm{cn}$ can be placed
arbitrarily by suitable choice of
$K_\mathrm{cn}$.  

The above statements on classical controllability %and observability
imply in particular that generically a single input %and a single output
($n_u=1$) can be used to stabilize the equilibrium
$(x_\mathrm{eq}(s),\mu_\mathrm{eq}(s))$ for any fixed $s$ (no matter
how many unstable directions $n_\mathrm{unst}(s)$ it has with $u=0$)
locally, by the scalar state feedback
%. For example, one may use observers $x_\mathrm{ob}$ and
%feedback control $u$ dynamically given by
%\cite{luenberger1971introduction}
\begin{align}
%\dot x_\mathrm{ob}&=f_xx_\mathrm{ob}+f_uu+K_\mathrm{ob}[g_x x_\mathrm{ob}-(y-y_\mathrm{eq})]\mbox{,}&
u&=K_\mathrm{cn}\cdot(x-x_\mathrm{eq})\mbox{.}
\label{classic:linfb}
\end{align}
The local exponential decay rate toward $x_\mathrm{eq}$ can be made
arbitrarily large with suitably chosen gains. % Since for $n_u=1$ linear
% controllability is a codimension-$1$ genericity condition,
% %(and similarly, for $n_y=1$ is linear observability),
% we have to expect that these may be violated along the equilibrium
% branch for isolated $s$.
In practice, stability becomes sensitive with
respect to the gains if one attempts to stabilize many unstable
degrees of freedom. Since the stabilizing gains depend on the partial
derivatives $f_x$, $f_u$ %, $g_x$
one needs good estimates for these. As we consider scenarios where we
do not have accurate estimates for partial derivatives, we will in
\Cref{sec:1dbranch} construct simple criteria for the gains in the single-input
single-output case $n_u=1$ and one degree of instability or less
($n_\mathrm{unst}\leq1$).

\subsubsection{Non-invasiveness}
\label{sec:non-inv}
As the concept of controllability recalled above is about linear
systems, when applying it to the vicinity of equilibria in nonlinear
systems one has to assume that the equilibrium location and the
partial derivatives $f_x$ and $f_u$ are known. Observe that
\eqref{classic:linfb} contains $x_\mathrm{eq}$ in its construction. In
practice, one constructs the feedback gains $K_\mathrm{cn}$ from estimates $\hat{f}_x$ and $\hat{f}_u$,
and inserts a \emph{reference value} $x_\mathrm{ref}$ into
\eqref{classic:linfb}:
\begin{align}\label{classic:statefb}
u&=K_\mathrm{cn}\cdot(x-x_\mathrm{ref})\mbox{.}
\end{align}
If the perturbations $\hat{f}_x-f_x$, $\hat{f}_u-f_u$ and
$x_\mathrm{ref}-x_\mathrm{eq}$ are sufficiently small then the
feedback gains $K_\mathrm{cn}$ constructed for
$\hat{f}_x$ and $\hat{f}_u$ are stabilizing the
equilibrium. That is, the controlled system
\eqref{eq:io:ode},\,\eqref{classic:statefb} will have
an equilibrium $x_\mathrm{cn}$ near $x_\mathrm{eq}$. Importantly,
\begin{align}
x_\mathrm{cn}=x_\mathrm{eq}+O(x_\mathrm{ref}-x_\mathrm{eq})\mbox{.}\label{eq:invasive}
\end{align}
That is, if  $x_\mathrm{ref}=x_\mathrm{eq}$ then small perturbation in the partial derivatives $\hat{f}_x-f_x$ and
$\hat{f}_u-f_u$ will not change the location of the equilibrium of the controlled system.

The application of feedback control to discover the precise equilibria
(or, more generally, steady states including periodic orbits)
motivates the concept of non-invasiveness.
When does it hold that $x_\mathrm{cn}=x_\mathrm{eq}$? Or, in other words, how to control the system such that the observed equilibrium coincides with the unknown equilibrium of the uncontrolled system?

\subsubsection{Non-invasiveness through iterative root finding}
\label{sec:gen:rootfinding}
If we assume controllability for all $s$ along the branch, and we
assume that the gain $K_\mathrm{cn}(s)\in \R^{1\times n_x}$ depends
smoothly on $s$ (for example, if it is constant), then the locally
stabilizing feedback control \eqref{classic:statefb},
$u(t)=K_\mathrm{cn}[x(t)-x_\mathrm{ref}]$,
for $x_\mathrm{ref}\approx x_\mathrm{eq}(s)$ induces
the input-output map
\begin{align}
  \label{gen:cbc:inoutmap} X_\infty: \R^{n_x}\times\R\ni(x_\mathrm{ref},\mu)\mapsto x_\mathrm{cn}=\lim_{t\to\infty}x(t)\in\R^{n_x}\mbox{.}
\end{align}
This map $X_\infty$ is evaluated at an argument $(x_0,\mu_0)$ near the
branch $(x_\mathrm{eq}(s),\mu_\mathrm{eq}(s))$ by setting the
reference value $x_\mathrm{ref}=x_0$ in \eqref{classic:statefb}, the
parameter $\mu=\mu_0$ in the system, \eqref{eq:io:ode}, waiting for
the transient dynamics of \eqref{eq:io:ode} to settle such that the
state $x(t)$ reaches a limit $x_\mathrm{cn}$. By construction, the
branch of equilibria $(x_\mathrm{eq}(s),\mu_\mathrm{eq}(s))$ are fixed
points of the map $X_\infty$, that is, they are solutions of the
equation
\begin{align}
  \label{gen:cbc:inout:fixedpoint}
  X_\infty(x,\mu)=x
\end{align}
for all $s\in[s_\mathrm{min},s_\mathrm{max}]$: when state $x(t)$
equals $x_\mathrm{ref}$ in \eqref{classic:statefb} then $u=0$ such that the
feedback control is non-invasive.

The papers
\cite{PhysRevLett.100.244101,doi:10.1177/1077546310384004,bureau2013experimental,barton2017control}
tracked branches of forced oscillations directly in mechanical
vibration experiments by applying linear feedback control (typically
proportional-plus-derivative feedback control
) to the forced nonlinear oscillators, and then used standard
numerical root-finding algorithms, such as simplified Newton
iterations, to find root curves of fixed point problem
\eqref{gen:cbc:inout:fixedpoint}. 

As implementing a stabilizing feedback loop is often the most
difficult part in physical experiments, permitting the experimenter to
choose gains with as few restrictions as possible is important
\cite{barton2017control}. The study \cite{bureau2013experimental}
performs a systematic investigation on how the gains can be chosen but
other experimental papers keep them constant along the branch. Using
maps $X_\infty$ where $u$ is based on adaptive feedback gains is
possible and improves robustness to uncertainty and time delays
\cite{li2021adaptive}. Its effect has been demonstrated on
ODEs similar to \eqref{eq:io:ode}.

A major obstacle for the application of standard root-finding
algorithms to solving \eqref{gen:cbc:inout:fixedpoint} is that
experimentally obtained data has high uncertainty such that precise
derivative information for the Jacobian of $X_\infty$ is not available
or expensive to obtain, especially in higher dimensions. This may lead
to slow and uncertain convergence, negating the main advantage of
classical numerical local root-finding and continuation algorithms
such as Newton iterations and pseudo-arclength continuation. Schilder
\emph{et al}.\ \cite{schilder2015experimental} developed and studied
modifications of these classical numerical methods, specifically to treat contamination with noise, which are now
available as \textsc{continex} toolbox in \textsc{coco} \cite{DS13}.

\subsection{Inherently non-invasive feedback control}
\label{sec:inh:noninv}
When the idea of finding unstable steady states (including periodic
orbits) by continuous-time feedback control was originally conceived,
the focus was on types of feedback that are inherently
non-invasive. These feedback laws introduce additional degrees of
freedom $x_\mathrm{wo}$. We will show that, in a suitable formulation,
these additional degrees of freedom act like integral components of a
proportional-plus-integral (PI) control, enforcing an additional
constraint, which is in this case $u=0$. With this formulation, it
becomes clear how to generalize existing types of inherently
non-invasive feedback to remove singularities at special points.
\subsubsection{Inherently non-invasive control --- washout filters}
\label{sec:SO}
For equilibria, the additional variables $x_\mathrm{wo}$ have been
introduced by \cite{Pyragas2004,WANG19951213} in the form of state
observers and were called washout filter in \cite{WANG19951213}
(hence the subscript $\mathrm{wo}$). Let us repeat the full non-linear
ODE such that the feature of non-invasiveness becomes clear (we keep
$s$ and, thus, $\mu$ fixed):
\begin{align}
  \label{gen:wo:ode}
  \dot x&=f(x,\mu,u)\mbox{,}&
  \dot x_\mathrm{wo}&=u\mbox{,\quad where $x_\mathrm{wo}\in\R^{n_u}$.}
\end{align}
The differential equation for $x_\mathrm{wo}$ implies that, whenever a linear feedback law of the general form
\begin{align}
  \label{gen:wo:u}
  u(t)=k_\mathrm{st}\cdot(x-x_\mathrm{ref})+k_\mathrm{wo}\cdot(x_\mathrm{wo}-x_\mathrm{wo,ref})\mbox{\ with arbitrary $k_\mathrm{st}\in\R^{n_u\times n_x}$ and $k_\mathrm{wo}\in\R^{n_u\times n_u}$}
\end{align}
is applied, every equilibrium of \eqref{gen:wo:ode},\,\eqref{gen:wo:u}
has $u=0$. Thus, equilibria of \eqref{gen:wo:ode},\,\eqref{gen:wo:u}
have a $x$-component that is also an equilibrium of the uncontrolled
system $\dot x=f(x,\mu,0)$ such that
\eqref{gen:wo:ode},\,\eqref{gen:wo:u} does not change the locations of
equilibria of the uncontrolled system. This justifies the notion of
inherent non-invasiveness for this type of feedback control. A simple
criterion, given in \Cref{thm:washout}, shows that a system is
controllable with washout filters whenever we have linear state
feedback controllability along the equilibrium branch
$(x_\mathrm{eq}(s),\mu_\mathrm{eq}(s))$, except when $f_x$ is
singular.
\begin{lemma}[Controllability for washout filters]\label{thm:washout}
  The linear time-invariant \textup{(}autonomous\textup{)} system
  \begin{align}
    \label{eq:thm:wo}
    \dot x&=Ax+Bu\mbox{,}& \dot x_\mathrm{wo}&=u&\mbox{with $A\in\R^{n_x\times n_x}$, $B\in\R^{n_x\times n_u}$,}
  \end{align}
  is controllable if and only if $\dot x=Ax+Bu$ is
  controllable and $A$ is regular.
\end{lemma}

\paragraph{Proof of \Cref{thm:washout}} The coefficients for the
state $x_\mathrm{ext}=(x,x_\mathrm{wo})$ and control $u$ in the right-hand side in
\eqref{eq:thm:wo} have the form ($I$ is the identity matrix)
\begin{align*}
  A_\mathrm{ext}&=
  \begin{bmatrix}
    A&0\\ 0&0
  \end{bmatrix}\mbox{,}&
  B_\mathrm{ext}&=
  \begin{bmatrix}
    B\\ I_{n_u\times n_u}
  \end{bmatrix}\mbox{.}
\end{align*}
Thus, the controllability matrix $R_\mathrm{ext}$ of the extended system has the form
\begin{align*}
  R_\mathrm{ext}=
  \begin{bmatrix}
    \begin{matrix}
      B\\  I_{n_u\times n_u}
    \end{matrix}&
    \begin{matrix}
      AB\\ 0
    \end{matrix}&\ldots&
    \begin{matrix}
      A^{n_x+n_u-1}B\\
      0
    \end{matrix}
  \end{bmatrix}\mbox{,}
\end{align*}
which has rank $n_x+n_u$ if and only if the matrix
$[AB,\ldots,A^{n_x}B]=A[B,AB,\ldots,A^{n_x-1}B]$ has rank $n_x$, where
$[B,AB,\ldots,A^{n_x-1}B]$ is the controllability matrix of
$\dot x=Ax+Bu$. \hfill$\square$

Formulation \eqref{gen:wo:ode},\,\eqref{gen:wo:u} is different from
the original papers \cite{AWC94,hassouneh2004washout,WANG19951213}. \js{Let us briefly explain that \eqref{gen:wo:ode},\,\eqref{gen:wo:u} encompasses the original washout filter formulations. The most general version for continuous-time washout filters introduced by Hassouneh \emph{et al}.\ \cite{hassouneh2004washout} is of the form (dropping nonlinear terms)
  \begin{align}\label{wo:hassouneh}
    \dot x&=Ax+Bu\mbox{,}&\dot z&=P(x-z)\mbox{,}& u=K(x-z)
  \end{align}
  with non-singular $P\in\R^{n_x\times n_x}$ and $z\in\R^{n_x}$. The
  authors showed that when the pair $(A,B)$ is stabilizable and $A$ is
  non-singular then one can find $P$, $K$ such that
  \eqref{wo:hassouneh} is asymptotically
  stable. System~\eqref{wo:hassouneh} is a special case of
  \eqref{eq:thm:wo} with a particular choice of gains if we assume
  that the gain $K$ in \eqref{wo:hassouneh} is a full-rank
  $n_u\times n_x$ matrix and that $n_u\leq n_x$. To see this, let us
  call $B_\mathrm{ext}=[B,0_{n_x\times(n_u-n_x)}]$ and choose a matrix $\tilde{K}\in \R^{(n_x-n_u)\times n_x}$ such that $[K^\tran,\tilde{K}^\tran]$ is non-singular. Then the variable
  \begin{align*}
    \begin{bmatrix}
      y\\\tilde{y}
    \end{bmatrix}=
    \begin{bmatrix}
      K\\\tilde{K}
    \end{bmatrix}
    P^{-1}z
  \end{align*}
  and $x$ satisfy
  \begin{align}\label{eq:wo_ext}
    \dot x&=A x+B_\mathrm{ext}
    \begin{bmatrix}
      u\\\tilde{u}
    \end{bmatrix}\mbox{,}&
    \begin{bmatrix}
      \dot y\\\Dot{\tilde{y}}
    \end{bmatrix}
&=
\begin{bmatrix}
  u\\
    \tilde{u}
  \end{bmatrix}
\end{align}
if $x$ and $z$ satisfy system
\eqref{wo:hassouneh}. System~\eqref{eq:wo_ext} is of the same form as
our formulation \eqref{eq:thm:wo} of the washout filter in
\Cref{thm:washout}, with $n_x-n_u$ components of the control input
unused. The form of system \eqref{wo:hassouneh} corresponds then to
the particular choice of gains
\begin{align*}
  k_\mathrm{st}&=
  \begin{bmatrix}
      K\\\tilde{K}
    \end{bmatrix}\mbox{,}&
    k_\mathrm{wo}&=\begin{bmatrix}
      K\\\tilde{K}
    \end{bmatrix}P\begin{bmatrix}
      K\\\tilde{K}
    \end{bmatrix}^{-1}
  \end{align*}
  in our formulation \eqref{gen:wo:ode}, \eqref{gen:wo:u}. Thus,
  \Cref{thm:washout} extends and simplifies the results of
  \cite{bazanella1997control,hassouneh2004washout}, as our result
  immediately implies all conclusions from controllability of linear systems with constant coefficients.} It
also separates the problem of controllability from the problem of
finding the control gains, which can then be designed using standard
linear feedback control theory.

\subsubsection{Inherently non-invasive feedback control through the parameter}
\label{sec:SP}
Siettos \emph{et al}.\ \cite{siettos2004coarse} showed that, if state
feedback control is applied through the bifurcation parameter $\mu$,
then an inherently non-invasive control can be constructed in the
context of a continuation of a branch of equilibria in $\mu$. The
feedback control in \cite{siettos2004coarse} considers $\mu$ as part
of the state, satisfying the
equation % It is known that the washout-filters as
% described in \eqref{eq:washout} are singular precisely at folds. To
% overcome this limitation, we adapt the additional differential
% equation to be a function of the system parameter $\mu$. It can be
% assumed to be
\begin{align}
  \label{siettos:mudot}
\dot{\mu}(t) &= 0 + u(t)\mbox{.}
\end{align}
In the setting of \cite{siettos2004coarse} this is the only point
where control input $u$ enters, such that this method considers a
scalar control input, $n_u=1$, for one-parameter equilibrium
branches. Thus, \eqref{eq:io:ode} has the form
\begin{align}
  \dot x(t)&=f(x(t),\mu(t))\mbox{.}\label{siettos:siso:ode}
\end{align}
Siettos \emph{et al}.\ observe that state feedback control of system
\eqref{siettos:mudot},\,\eqref{siettos:siso:ode} of the form
\begin{align}\label{siettos:u}
  u(t)=K_{\mathrm{st},x}(x(t)-x_\mathrm{ref})+K_{\mathrm{st},\mu}(\mu(t)-\mu_\mathrm{ref})
\end{align}
is always non-invasive in the sense that equilibria of the controlled
system must satisfy $u=0$, such that they will lie on the intersection
of the equilibrium curve
$(x_\mathrm{eq}(\cdot),\mu_\mathrm{eq}(\cdot))$ with the codimension-$1$ hyperplane $u=0$ in
the $(x,\mu)$-space.
Siettos \emph{et al}.\ \cite{siettos2004coarse} demonstrated
continuation through their feedback control for a kinetic Monte-Carlo
simulation.  A control law similar to \eqref{siettos:u} was proposed
for Poincare maps of periodic orbits in \cite{
  misra2008event}. Physical experiments on vibrations of nonlinear
mechanical oscillators were performed by \cite{BartonSieber2013} on a
further simplification of
\eqref{siettos:mudot},\,\eqref{siettos:siso:ode},\,\eqref{siettos:u} by applying the
feedback control law
\begin{align}
  \mu(t)=\mu_\mathrm{ref}+K_\mathrm{st}\cdot(x(t)-x_\mathrm{ref})\mbox{,}\label{eq:muBS}
\end{align}
which corresponds to choosing
$|K_{\mathrm{st},\mu}|,|K_{\mathrm{st},x}|\gg1$ with fixed ratio
$K_\mathrm{st}=K_{\mathrm{st},x}/K_{\mathrm{st},\mu}$) in
\eqref{siettos:u}. This law is also non-invasive whenever it is
stabilizing.
\begin{lemma}[Controllability for parameter control]\label{thm:SP}\ \\
  An equilibrium $(x_\mathrm{eq}(s),\mu_\mathrm{eq}(s))$ of
  \eqref{siettos:mudot}, \eqref{siettos:siso:ode} is
  controllable if and only if the pair of partial derivatives
  $(f_x(s),f_\mu(s))$ in $(x_\mathrm{eq}(s),\mu_\mathrm{eq}(s))$ is
  controllable.
\end{lemma}
\paragraph{Proof of \Cref{thm:SP}} The coefficients for the
state $x_\mathrm{ext}=(x,\mu)$ and control $u$ in the right-hand side in
\eqref{siettos:mudot},\,\eqref{siettos:siso:ode} have the form
\begin{align*}
  A_\mathrm{ext}&=
  \begin{bmatrix}
    f_x&f_\mu\\ 0&0
  \end{bmatrix}\mbox{,}&
  B_\mathrm{ext}&=
  \begin{bmatrix}
    0\\ 1
  \end{bmatrix}\mbox{.}
\end{align*}
Thus, the controllability matrix $R_\mathrm{ext}$ of the extended system has the form
\begin{align*}
  R_\mathrm{ext}=
  \begin{bmatrix}
    \begin{matrix}
      0\\ 1
    \end{matrix}&
    \begin{matrix}
      f_\mu\\ 0
    \end{matrix}&
        \begin{matrix}
      f_xf_\mu\\ 0
    \end{matrix}&\ldots&
    \begin{matrix}
      f_x^{n_x-1}f_\mu\\
      0
    \end{matrix}
  \end{bmatrix}\mbox{.}
\end{align*}
This matrix has rank $n_x+1$ if and only if the
matrix $[f_\mu,\ldots,f_x^{n_x-1}f_\mu]$ has rank $n_x$, which
is the controllability matrix of the system
$\dot x=f_xx+f_\mu u$. \hfill$\square$

Thus, the controllability condition for \eqref{siettos:u} is the same
as for the simplified control law \eqref{eq:muBS}.

\subsubsection{Inherently non-invasive control -- zero-in-equilibrium feedback}
\label{sec:combi}
Both methods for inherently non-invasive control of equilibria are
expected to fail at isolated points along the one-parameter
equilibrium branch $(x_\mathrm{eq}(s),\mu_\mathrm{eq}(s))$ for 
\begin{align}
  \label{combi:ode}
  \dot x=f(x,\mu,u)\mbox{.}
\end{align}
For controllability through the inputs $\mu$ or $u$ on their
own, the relevant controllability matrices are (dropping argument $s$)
\begin{align}
  \label{combi:Ru}
  R_u&=[f_u,f_xf_u,\ldots,f_x^{n_x-1}f_u]\in\R^{n_x\times(n_x\cdot n_u)}\mbox{,}\\
  \label{combi:Rmu}
  R_\mu&=[f_\mu,f_xf_\mu,\ldots,f_x^{n_x-1}f_\mu]\in\R^{n_x\times n_x}\mbox{.}  
\end{align}
\begin{itemize}
\item \textbf{Failure at singular $f_x(s)$:} Washout filter based
  control with input $u$ for \eqref{combi:ode} % , such as washout
  % filters and time-delayed feedback control,
  do not stabilize the
  equilibrium near fold bifurcations. This failure occurs even if the
  pair $(f_x(s),f_u(s))$ is controllable in the parameter $s$ of the
  fold bifurcation, that is, when $R_u(s)$ has rank $n_x$.
\item \textbf{Failure at singular $R_\mu(s)$:} For one-parameter
  branches the control input for parameter control is naturally
  one-dimensional ($n_u=1$) such that a loss of controllability
  (singularity of $R_\mu(s)$) is a codimension-$1$ event, expected to
  occur at isolated points along the branch
  $(x_\mathrm{eq}(s),\mu_\mathrm{eq}(s))$.
\end{itemize}
In the particle flow model example in \Cref{sec:scenario} the time
scale of the response of the particle flow to parameter changes (in
our case the position $\mu$ of the obstacle) is also too slow, such
that non-invasive feedback control, when applied purely through the bifurcation parameter,
creates a large uncertainty in the resulting steady states
$x_\mathrm{eq}$, observed as limits $\lim_{t\to\infty}x(t)$ after
transients have settled. See also \Cref{sec:control-based-ped} where we briefly discuss control through the bifurcation parameter.

Hence, we extend feedback control with washout filters by using the
bifurcation parameter $\mu$ as the integral component, in the same way
as in the feedback control through the parameter. This permits us to
design an inherently non-invasive control that is only singular at
events of codimension larger or equal than $2$. % The approach is especially
% suitable for computational experiments such as our particle flow model
% because additional control inputs $u$ can be added with little
% cost.
Let us consider the feedback control scheme
\begin{align}
  \label{combi:ode+control}
  \dot x&=f(x,\mu,au)\mbox{,}&
  \dot \mu&=u\mbox{\quad with $\dim u=1$ (thus, $n_u=1$), and $a\in\R$,}
\end{align}
such that $u$ is an input in $f$ in addition to the also
adjustable $\mu$.  The scalar $a$ acts as an additional weight on the
control gains in $f$ compared to $\dot\mu$, which we are free to
choose suitably. We observe that this feedback is also inherently
non-invasive: every equilibrium of \eqref{combi:ode+control} is also
an equilibrium of the uncontrolled system \eqref{combi:ode} with
$u=0$. We refer to \eqref{combi:ode+control} as \emph{zero-in-equilibrium feedback control}.
\begin{lemma}[Controllability for %non-invasive control using extended washout filter
zero-in-equilibrium feedback control]\label{thm:combi}
  \ \\ The equilibrium $(x_\mathrm{eq},\mu_\mathrm{eq})$ of
  \eqref{combi:ode+control} is controllable if the matrix
  $a f_xR_u+R_\mu$ is regular, where $R_u$ and $R_\mu$ are the controllability matrices defined in \eqref{combi:Ru} and \eqref{combi:Rmu}.
\end{lemma}
\paragraph{Proof of \Cref{thm:combi}} The coefficients of the linearization of system \eqref{combi:ode+control} are
\begin{align*}
  A_\mathrm{ext}&=
  \begin{bmatrix}
    f_x&f_\mu\\
    0&0
  \end{bmatrix}\mbox{,}&
  B_\mathrm{ext}&=
  \begin{bmatrix}
    a f_u\\
    1
  \end{bmatrix}\mbox{.}
\end{align*}
Thus, the controllability matrix for the linearized system is
\begin{align*}
  R_\mathrm{ext}&=
  \begin{bmatrix}
    \begin{matrix}
    a f_u\\ 1
  \end{matrix}&
  \begin{matrix}
    a f_xf_u+f_\mu\\ 0
  \end{matrix}&
  \begin{matrix}
    a f_x^2f_u+f_xf_\mu\\ 0
  \end{matrix}& \ldots
    \begin{matrix}
    a f_x^{n_x}f_u+f_x^{n_x-1}f_\mu\\ 0
  \end{matrix}&
\end{bmatrix}\mbox{,}
\end{align*}
which has rank $n_x+1$ if and only if $a f_x\cdot R_u+R_\mu$ is regular.\hfill$\square$

When constructing non-invasive control, one faces the question when
one can find a scalar $a$ for which the equilibrium
$(x_\mathrm{eq},\mu_\mathrm{eq})$ is controllable. We
may phrase the answer given by \Cref{thm:combi} in terms of regularity
of the matrix pair $(R_\mu,f_xR_u)$ of $\R^{n_x\times n_x}$ matrices.
A pair $(A_0,A_1)$ of $\R^{n_x\times n_x}$ matrices is called regular, if
the polynomial $\lambda\mapsto\det(A_0+\lambda A_1)$ is not identically
zero. The condition on the coefficients in the two matrices $A_0$ and
$A_1$ for the pair to be singular (i.e., not regular) imposes $n_x+1$
constraints: all coefficients of the characteristic polynomial
$\lambda\mapsto\det(A_0+\lambda A_1)$ have to be zero. Thus, violations of
matrix pair regularity are codimension $n_x+1$
events. \Cref{thm:combi:controllability} summarizes the consequences
of linear feedback controllability for control law
\eqref{combi:ode+control}.
\begin{corollary}[Existence of stabilizing gains]\label{thm:combi:controllability}
  Let  $\gamma>0$ be
  arbitrary. If the matrix pair $(R_\mu,f_xR_u)$ is regular
  for the linearization of $\dot x=f(x,\mu,0)$ in equilibrium
  $(x_\mathrm{eq},\mu_\mathrm{eq})$, then there exist control
  gains $K_{\mathrm{st},x}\in\R^{1\times n_x}$,
  $K_{\mathrm{st},\mu}\in\R$, $a\in\R$, such that the  controlled
  system $\dot x=f(x,\mu,au)$, $\dot\mu=u$ with
   \begin{align}\label{combi:u}
     u=K_{\mathrm{st},x}(x_\mathrm{ref}-x)+K_{\mathrm{st},\mu}(\mu_\mathrm{ref}-\mu)
   \end{align}
   has for all
   $(x_\mathrm{ref},\mu_\mathrm{ref})\approx(x_\mathrm{eq},\mu_\mathrm{eq})$
   a stable equilibrium
   $(x_\mathrm{cn},\mu_\mathrm{cn})\approx(x_\mathrm{eq},\mu_\mathrm{eq})$, which is reached 
   with local exponential decay rate greater than $\gamma$.
\end{corollary}

\paragraph{Proof of \Cref{thm:combi:controllability}}
We choose the weight $a$ such that $af_xR_u+R_u$ is regular, which is possible
by the regularity of the matrix pair. The resulting controllability of
the linearization in $(x_\mathrm{eq},\mu_\mathrm{eq})$ permits us to
choose gains $(K_{\mathrm{st},x},K_{\mathrm{st},\mu})$ such that the
linearization of system \eqref{combi:ode+control} in
$(x_\mathrm{eq},\mu_\mathrm{eq})$ with feedback control $u$ given in
\eqref{combi:u} and
$(x_\mathrm{ref},\mu_\mathrm{ref})=(x_\mathrm{eq},\mu_\mathrm{eq})$
has a spectrum where all eigenvalues have real part less than
$-\gamma$. Continuity then ensures that the decay rates and the
equilibrium persist for $(x_\mathrm{ref},\mu_\mathrm{ref})$ near
$(x_\mathrm{eq},\mu_\mathrm{eq})$.\hfill$\square$

In contrast to the washout filters \eqref{gen:wo:ode} or control through the
bifurcation parameter \eqref{siettos:siso:ode},\,\eqref{siettos:u},
for which controllability conditions fail at events of codimension
$1$, control law \eqref{combi:ode+control}, using the bifurcation
parameter as observer and an additional control input $u$, fails only
at events of codimension $n_x+1$. 

Adding a real-time feedback control input $u$ imposes a cost in
physical experiments. However, in computational experiments such as
our particle flow model for a pedestrian evacuation scenario, an
additional input has negligible cost and can be constructed to make
choosing stabilizing gains $(a,K_{\mathrm{st},x},K_{\mathrm{st},\mu})$
as easy as possible.

\subsection{Branches with single slow dimension}
\label{sec:1dbranch}
The statements in \Cref{sec:inh:noninv} are
concerned with non-invasive controllability of systems with arbitrary
state dimension $n_x$. In particular, they permit an arbitrary number
$n_\mathrm{unst}$ of unstable dimensions for the equilibrium
$(x_\mathrm{eq},\mu_\mathrm{eq})$. General control theory also gives
explicit procedures to construct gains (such as
$(K_{\mathrm{st},x},K_{\mathrm{st},\mu})$ in \eqref{siettos:u} or
\eqref{combi:u}) resulting in arbitrary decay rates toward the controlled
equilibrium. \js{However, these procedures rely on precise knowledge of the partial derivatives $f_x$ and $f_u$, and the results are sensitive to errors in estimating these derivatives. For this reason we now discuss the common scenario that $n_\mathrm{unst}\leq 2$ and that $f_x$, $f_\mu$,
$f_u$ or $g_x$ are difficult or computationally expensive to
approximate. For these cases we can state simple inequality constraints on the scalar gains that ensure stabilization.}

In particular we hypothesized that two fold
bifurcations and a branch of unstable equilibria 
cause the bistability in our multi-particle model for pedestrians, shown in \Cref{fig:intro}.
This implicitly includes the hypothesis that the model behaves essentially
as a system of ODEs close to a branch of equilibria where the number
of dimensions changing stability equals $1$, and where all other
eigenvalues in the equilibria have uniformly negative real part. The
fluctuations observed around a stable stationary particle flow are
then treated as perturbations generating uncertainty, similar to a
physical experiment.

For this case of an essentially one-dimensional ODE the criteria for
gains in the non-invasive control schemes discussed in
\Cref{sec:gendesign} and \Cref{sec:combi} to achieve at least
stabilization can be simplified and made explicit, which we will
discuss in this section.

We consider the system with output
\begin{align}
  \label{1d:ode+out} \dot x&=f(x,\mu,u)\mbox{,}&y&=g(x)\mbox{\quad with $n_y=\dim y=n_u=\dim u=1$}
\end{align}
for the case where along the branch of equilibria
($x_\mathrm{eq}(s),\mu_\mathrm{eq}(s))$ the Jacobian $f_x(s)$ has only
one direction $v_\mathrm{c}(s)$ in which the growth rate
$\lambda_\mathrm{c}(s)$ is of order $1$, while all other directions
are strongly stable with time scale difference of order
$\epsilon$. More precisely:
\begin{assumption}[Time scale difference and single slow (center) direction]
  \label{ass:1d}
  We assume that the Jacobian $f_x(s)$ in the equilibrium
  $(x_\mathrm{eq}(s),\mu_\mathrm{eq}(s))$ of \eqref{1d:ode+out} with
  $u=0$ has a single simple eigenvalue $\lambda_\mathrm{c}(s)\in\R$
  with right and left eigenvectors $v_\mathrm{c}(s)$,
  $w_\mathrm{c}(s)\in\R^{n\times1}$ with modulus of order $O(1)$,
  while all others are stable with time scale ratio
  $\epsilon$\textup{:}
\begin{align}
  \label{ass:center} &f_x v_\mathrm{c}=\lambda_\mathrm{c}v_\mathrm{c}\mbox{,\quad}
  w_\mathrm{c}^\tran{}f_x=\lambda_\mathrm{c}w_\mathrm{c}^\tran{}
  \mbox{\quad with scaling\ } 1=w_\mathrm{c}^\tran{}v_\mathrm{c}\mbox{,\ }|\lambda_\mathrm{c}|=O(1)\mbox{, and}\\
  %\label{ass:scaling} 1&=w_\mathrm{c}^\tran{}v_\mathrm{c}\mbox{,}& |\lambda_\mathrm{c}|&=O(1)\mbox{,\quad while, for}\\
  \label{ass:stab}&\re\left[\spec f_x\vert_{\ker w_\mathrm{c}^\tran{}}\right]<-c_\mathrm{spec}/\epsilon\mbox{\quad and\quad}
  \left\|\left[f_x\vert_{\ker w_\mathrm{c}^\tran{}}\right]^{-1}\right\|\leq \epsilon c_\mathrm{st}\mbox{}
\end{align}
for some positive constants $c_\mathrm{spec}$ and $c_\mathrm{st}$ of order $1$ and $\epsilon\ll1$.
\end{assumption}
All variables in \Cref{ass:1d}, $\lambda_\mathrm{c}$, $v_\mathrm{c}$,
$w_\mathrm{c}$ and $f_x\vert_{\ker w_\mathrm{c}^\tran{}}$, depend on
$s$ but $c_\mathrm{st}$, $c_\mathrm{spec}$ and $\epsilon$ are independent of $s$.  In \eqref{ass:stab}
%\begin{align*}
  $f_x(s)\vert_{\ker w_\mathrm{c}^\tran{}(s)}$ % =:f_{x,\mathrm{st}}(s) 
%\end{align*}
is the stable part of the Jacobian $f_x(s)$.
\begin{assumption}[Partial linear observability of slow direction]
\label{ass:1d:obs:cont}  
We  assume that for all $s\in[s_\mathrm{min},s_\mathrm{max}]$
\begin{align}
  \label{1d:linobservability}
\R&\ni\ \, g_x(s) v_\mathrm{c}(s)\neq 0   &&\mbox{\textup{(}partial observability\textup{).}}%\\
%  \label{eq:lincontrollability}
%\R&\ni w_\mathrm{c}^\tran{}(s)f_u(s)\neq 0 %\mbox{, and\quad}
%\R\ni w_\mathrm{c}^\tran{}(s)f_\mu(s)\neq 0
%&&\mbox{\textup{(}partial controllability\textup{)}.}
%\mbox{}
\end{align}
\end{assumption}
Thus, the scalar quantity $g_x(s) v_\mathrm{c}(s)$ never changes sign
along the branch and we may scale the eigenvector $v_\mathrm{c}(s)$
such that
\begin{align}
  \label{1d:scale:obs} 1=g_x(s) v_\mathrm{c}(s)\mbox{\quad for all $s\in[s_\mathrm{min},s_\mathrm{max}]$.}
\end{align}
Assumption \eqref{1d:linobservability} is a weaker genericity
assumption than full linear observability, as we only want to observe
the slow direction $v_\mathrm{c}$ through output $y=g(x)$.  In
contrast,  assumptions such as
\begin{align}
  \label{1d:lincontrollability:u}
\R&\ni w_\mathrm{c}^\tran{}(s)f_u(s)\neq 0 &&\mbox{(partial controllability through $u$), or}\\
  \label{1d:lincontrollability:mu}
\R&\ni w_\mathrm{c}^\tran{}(s)f_\mu(s)\neq 0
&&\mbox{(partial controllability through $\mu$),}
\end{align}
are not necessarily weaker genericity assumptions than the respective
full controllability, as we do not want to rely on the coupling
from stable directions in $\ker w_\mathrm{c}^\tran{}$ for
stabilizing the equilibrium in the $v_\mathrm{c}$ direction.
So, we are not making controllability assumptions at this stage but will consider them later for each particular type of non-invasive control laws.

In the $\epsilon$-vicinity of an equilibrium
$(x_\mathrm{eq}(s),\mu_\mathrm{eq}(s))$ we split the deviation of the
state $x(t)$ from its equilibrium into its slow and its stable parts, %  and
% name the deviations of possibly varying $\mu$ and $y$ from their
% equilibrium value
\begin{align}
  \label{ass:xsplit} x_\mathrm{dev}(t)&:=x(t)-x_\mathrm{eq}(s)=v_\mathrm{c}(s)x_\mathrm{c}(t)+V_\mathrm{stb}(s)x_\mathrm{stb}(t)\mbox{.}
\end{align}
(dropping the argument $s$ from the deviation) where the rows of
$V_\mathrm{stb}(s)\in\R^{n\times(n-1)}$ span
$\ker w_\mathrm{c}^\tran{}(s)$ (i.e.,
$0=w_\mathrm{c}^\tran{}(s)V_\mathrm{stb}(s)$,
$I=V_\mathrm{stb}^\tran{}(s)V_\mathrm{stb}(s)\in\R^{(n-1)\times(n-1)}$)
and $x_\mathrm{c}(t)\in\R$, $x_\mathrm{stb}\in\R^{n-1}$ with
$|x_\mathrm{c}|$, $\|x_\mathrm{stb}\|=O(\epsilon)$.

The different non-invasive control laws, discussed in
\Cref{sec:inh:noninv}, applied to system \eqref{1d:ode+out} with
scalar output and single slow dimension (\Cref{ass:1d}) result in the
expressions and conditions for control gains, discussed in the
following paragraphs.

\subsubsection{Non-invasive control based on washout filters}
\label{sec:1d:SO}
 The general principle was discussed in \Cref{sec:SO}. For partial controllability through input $u$, \eqref{1d:lincontrollability:u}
($w_\mathrm{c}^\tran{}f_u\neq0$), we construct the washout filter and feedback as
\begin{align}  
  \label{1d:SO}
  \dot y_\mathrm{wo}&=u\mbox{,}&u&=K_\mathrm{st}y+K_\mathrm{wo}y_\mathrm{wo}
\end{align}
(setting $y_\mathrm{ref}=y_\mathrm{wo,ref}=0$ without loss of generality). If
$K_\mathrm{st}$ and $K_\mathrm{wo}$ are of order $O(1)$,
there exists a two-dimensional invariant slow manifold. The slow
coordinates $y$ and $y_\mathrm{wo}$ satisfy to first order the
equation
\begin{align*}
    \dot y&=\lambda_\mathrm{c}[y-y_\mathrm{eq}]+
    w_\mathrm{c}^\tran{}f_u\left[K_\mathrm{st}y+K_\mathrm{wo}y_\mathrm{wo}\right]+O(\|(y-y_\mathrm{eq},u)\|^2)\mbox{,}\\
    \dot y_\mathrm{wo}&=K_\mathrm{st}y+K_\mathrm{wo}y_\mathrm{wo}\mbox{,}
\end{align*}
which has the Jacobian in the equilibrium
$y=y_\mathrm{eq}$, $u=0$
\begin{align*}
  A_\mathrm{wo}=
  \begin{pmatrix*}[c]
    \lambda_\mathrm{c}+w_\mathrm{c}^\tran{}f_uK_\mathrm{st} & w_\mathrm{c}^\tran{}f_uK_\mathrm{wo} \\
    K_\mathrm{st} & K_\mathrm{wo}      
  \end{pmatrix*}\mbox{.}      
\end{align*}
Thus, the equilibrium is stable, if $A_\mathrm{wo}$ satisfies
$\tr A_\mathrm{wo}=\lambda_\mathrm{c}+K_\mathrm{wo}+w_\mathrm{c}^\tran{}f_uK_\mathrm{st}<0$, $\det A_\mathrm{wo}=\lambda_\mathrm{c}K_\mathrm{wo}>0$, which are
equivalent to
\begin{align}\label{washout:unstable}
  \lambda_\mathrm{c}+K_\mathrm{wo}&<(-w_\mathrm{c}^\tran{}f_u)K_\mathrm{st}\mbox{,} &K_\mathrm{wo}&>0\mbox{,}
  &&\mbox{if $x_\mathrm{eq}$ is unstable ($\lambda_\mathrm{c}>0$),}\\
  \label{washout:stable}
  \lambda_\mathrm{c}+K_\mathrm{wo}&<(-w_\mathrm{c}^\tran{}f_u)K_\mathrm{st}\mbox{,} &K_\mathrm{wo}&<0\mbox{,}
  &&\mbox{if $x_\mathrm{eq}$ is stable$\phantom{\mathrm{un}}$ ($\lambda_\mathrm{c}<0$).}
\end{align}
Thus, the sign of $K_\mathrm{wo}$ has to be chosen depending on the
stability of the branch (with arbitrary modulus, e.g.,
$K_\mathrm{wo}=\pm1$). After choosing the sign of $K_\mathrm{st}$ %has to be chosen
suitably, then the modulus of $K_\mathrm{st}$ has to be chosen
sufficiently large. 
\js{Consequently, a sufficient non-degeneracy condition for the existence of stabilizing gains $(K_{\mathrm{st}},K_{\mathrm{wo}})$ is that
\begin{equation}
  \label{1d:cond:wo}
\lambda_\mathrm{c}\neq 0.
\end{equation}}

It is also clear that $A_\mathrm{wo}$ is singular
if $\lambda_\mathrm{c}=0$ such that the 
non-invasive feedback control based on washout filter fails for all possible gains at fold bifurcations (when $\lambda_\mathrm{c}=0$).

\subsubsection{Non-invasive control through the bifurcation parameter}
\label{sec:1d:SP}
 The general principle was discussed in \Cref{sec:SP}. We do not assume the presence of an input $u$ in
the right-hand side (thus, $\dot x=f(x,\mu,0)$, $y=g(x )$), but require partial controllability through the bifurcation parameter $\mu$,
\eqref{1d:lincontrollability:mu}
($w_\mathrm{c}^\tran{}f_\mu\neq0$), and set
\begin{align}
    \label{1d:SP}
  \dot\mu&=u=K_{\mathrm{st},y}[y-y_\mathrm{ref}]+K_{\mathrm{st},\mu}[\mu-\mu_\mathrm{ref}]\mbox{.}
\end{align}
\js{If $K_{\mathrm{st},y}$ and $K_{\mathrm{st},\mu}$ are of order $O(1)$, and
$y_\mathrm{ref}$ and $\mu_\mathrm{ref}$ are near $y_\mathrm{eq}$ and
$\mu_\mathrm{eq}$,
there exists a two-dimensional invariant slow manifold. The slow
coordinates $y$ and $\mu$ satisfy to first order the
equation}
\begin{align*}
    \dot y&=\lambda_\mathrm{c}[y-y_\mathrm{eq}]+
    w_\mathrm{c}^\tran{}f_\mu [\mu-\mu_\mathrm{eq}]
    +O(\|(y-y_\mathrm{eq},\mu-\mu_\mathrm{eq})\|^2)\mbox{,}\\
    \dot \mu&=K_{\mathrm{st},y}[y-y_\mathrm{ref}]+K_{\mathrm{st},\mu}[\mu-\mu_\mathrm{ref}]\mbox{,}
\end{align*}
which has the Jacobian in the equilibrium
$(x_\mathrm{eq},\mu_\mathrm{eq})$
\begin{align*}
  A_\mu=
  \begin{pmatrix*}[c]
    \lambda_\mathrm{c}& w_\mathrm{c}^\tran{}f_\mu \\
    K_{\mathrm{st},y} & K_{\mathrm{st},\mu}
  \end{pmatrix*}\mbox{.}      
\end{align*}
Thus, criteria for stabilizing gains are that
\begin{align}\label{1d:SP:gaincrit1}
 K_{\mathrm{st},\mu}&<-\lambda_\mathrm{c}\mbox{,}&
\lambda_\mathrm{c} K_{\mathrm{st},\mu}-(w_\mathrm{c}^\tran{}f_\mu) K_{\mathrm{st},y}>0\mbox{.}
\end{align}
\js{
Consequently, a sufficient non-degeneracy condition for the existence of stabilizing gains $(K_{\mathrm{st},y},K_{\mathrm{st},\mu})$ is that
\begin{align}
  \label{1d:combi:gaincrit}
  w_\mathrm{c}^\tran{}f_\mu&\neq0
  \mbox{\qquad or\qquad}\lambda_\mathrm{c}<0\mbox{,}
\end{align}
}
If the equilibrium is part of a branch
$(x_\mathrm{eq}(s),\mu_\mathrm{eq}(s))$, then the ratio
$(\lambda_\mathrm{c})/(w_\mathrm{c}^\tran{}f_\mu)$ present in the
second condition in \eqref{1d:SP:gaincrit1} has a geometric
interpretation under one additional assumption: differentiating the
identity for equilibria, $f(x_\mathrm{eq}(s),\mu_\mathrm{eq}(s))=0$,
with respect to $s$ and projecting the resulting linear relation
between $\partial_s x_\mathrm{eq}$ and $\partial_s \mu_\mathrm{eq}$ by
$w_\mathrm{c}^\tran{}$, we obtain the linear relation
\begin{align}\label{1d:eqcenter:x}
  \lambda_\mathrm{c}w_\mathrm{c}^\tran{}\partial_sx_\mathrm{eq}+(w_\mathrm{c}^\tran{}f_\mu)\partial_s\mu_\mathrm{eq}=0\mbox{.}
\end{align}
If we assume in addition that the spectral stable projection of $f_\mu$ is not large, that is,
\begin{align}
  \label{1d:outcentral}
  [I-v_\mathrm{c}w_\mathrm{c}^\tran{}]f_\mu=O(1)\mbox{,}
\end{align}
($[I-v_\mathrm{c}w_\mathrm{c}^\tran{}]$ is the spectral projection
for $f_x$ onto $\ker w_\mathrm{c}^\tran{}$) then, by \Cref{ass:1d},
stability of $f_x\vert_{\ker w_\mathrm{c}^\tran{}}$ with timescale $1/\epsilon$, \eqref{ass:stab},
$[I-v_\mathrm{c}w_\mathrm{c}^\tran{}]\partial_sx_\mathrm{eq}=[f_x\vert_{\ker
  w_\mathrm{c}^\tran{}}]^{-1}[I-v_\mathrm{c}w_\mathrm{c}^\tran{}]f_\mu=O(\epsilon)\ll1$. Consequently,
\begin{align}\label{1d:ycenter}
\partial_s y_\mathrm{eq}=g_x\partial_s
x_\mathrm{eq}=g_xv_\mathrm{c}w_\mathrm{c}^\tran{}\partial_s
x_\mathrm{eq}+O(\epsilon)=w_\mathrm{c}^\tran{}\partial_s
x_\mathrm{eq}+O(\epsilon)\mbox{.}
\end{align}
Inserting $\partial_s y_\mathrm{eq}$ for $w_\mathrm{c}^\tran{}\partial_s
x_\mathrm{eq}$ in \eqref{1d:eqcenter:x}, results in the relation
\begin{align}\label{1d:eqcenter}
  \lambda_\mathrm{c}\partial_sy_\mathrm{eq}+(w_\mathrm{c}^\tran{}f_\mu)\partial_s\mu_\mathrm{eq}=O(\epsilon)\mbox{.}
\end{align}
Thus, the second condition on the gains to be stabilizing can be
phrased in terms of the tangent of the equilibrium curve in the
$(y,\mu)$-plane, $(y_\mathrm{eq}(s),\mu_\mathrm{eq}(s))$. The two vectors
$(\lambda_\mathrm{c},w_\mathrm{c}^\tran{}f_\mu)$ and
$(\partial_sy_\mathrm{eq},\partial_s\mu_\mathrm{eq})$ are both
non-zero and approximately orthogonal to each other along the
equilibrium branch. Along a
stable part of the equilibrium branch (where $\lambda_\mathrm{c}<0$
and the signs of $\partial_sy_\mathrm{eq}$ and
$\partial_s\mu_\mathrm{eq}$ can be established with zero control
gains), we may establish a sign
$\sigma=\pm1$ (independent of $s$) such that there exists
a $p(s)>0$ with
$(\lambda_\mathrm{c},w_\mathrm{c}^\tran{}f_\mu)=p\sigma(\partial_s\mu_\mathrm{eq},-\partial_sy_\mathrm{eq})+O(\epsilon)$
for all $s\in[s_\mathrm{min},s_\mathrm{max}]$. 
Thus, overall the criteria for the gains are
\begin{align}\label{1d:SP:gaincrit1:final}
 K_{\mathrm{st},\mu}&<-\lambda_\mathrm{c}\mbox{,}&
\sigma\left[\partial_s\mu_\mathrm{eq}K_{\mathrm{st},\mu}+\partial_sy_\mathrm{eq}K_{\mathrm{st},y}\right]+O(\epsilon)&>0\mbox{.}
\end{align}
In other words, the gains $(K_{\mathrm{st},y},K_{\mathrm{st},\mu})$
need to be sufficiently large in modulus and the line in the $(y,\mu)$
plane defined by
$0=u=K_{\mathrm{st},y}[y-y_\mathrm{ref}]+K_{\mathrm{st},\mu}[\mu-\mu_\mathrm{ref}]$
must intersect the equilibrium curve
$(y_\mathrm{eq}(s),\mu_\mathrm{eq}(s))$ at a non-zero angle %\js{, see also \cref{fig:geometric_intuition}}
. The orientation is determined by $\sigma$, such that we call the sign $\sigma$ the \emph{input orientation}.
%\begin{figure}[t]
%\centering
%\label{fig:geometric_intuition}
%%\includegraphics[scale=0.1]{Figures/geom_int.jpg}
%\includegraphics[scale=1]{Figures/geometric_intuition}
%\caption{\js{\jsq{I would have thought that this is not necessary if fig~\ref{fig:flow_manipulation_js} is sufficiently clear.} Illustration of the geometric intuition behind the control through the bifurcation parameter. The control gains $(K_{\mathrm{st},y},K_{\mathrm{st},\mu})$
%need to be chosen such that the line in the $(y,\mu)$
%plane defined by
%$K_{\mathrm{st},y}[y-y_\mathrm{ref}]+K_{\mathrm{st},\mu}[\mu-\mu_\mathrm{ref}]=0$ intersects the equilibrium curve
%$(y_\mathrm{eq}(s),\mu_\mathrm{eq}(s))$ at a non-zero angle.}%\jsq{I think that your middle panel in \cref{fig:flow_manipulation} can be modified to do the illustration.}
%}
%\end{figure}

The additional condition \eqref{1d:outcentral} is best understood by
its primary consequence \eqref{1d:ycenter}. The tangent to the
equilibrium curve $(x_\mathrm{eq}(s),\mu_\mathrm{eq}(s))$ should be
mostly tangential to the $( y,\mu)$-plane. Thus, changes in equilibrium
location and slow dynamics should be approximately aligned. This
condition is known to be satisfied at a fold bifurcation in $\mu$.

\subsubsection{Non-invasive control %through extended washout filter
based on zero-in-equilibrium feedback}
\label{sec:1d:combi}
The general principle was discussed in  \Cref{sec:combi}. The difference to \Cref{sec:1d:SP} is that the
input $u$ in the right-hand side is present, such that $\dot x=f(x,\mu,au)$ with non-zero $a$, $y=g(x(t))$. We set, identically to  Eq.~\eqref{1d:SP}, 
\begin{align}
  \label{1d:combi}
  \dot\mu&=u\mbox{,}&u&=K_{\mathrm{st},y}[y-y_\mathrm{ref}]+K_{\mathrm{st},\mu}[\mu-\mu_\mathrm{ref}]\mbox{.}
\end{align}
%and scale the input $u$ in the right-hand side, such that $\dot x=f(x,\mu,au)$.
\js{If $a$, $K_{\mathrm{st},y}$ and $K_{\mathrm{st},\mu}$ are of order $O(1)$, and
$y_\mathrm{ref}$ and $\mu_\mathrm{ref}$ are near $y_\mathrm{eq}$ and $\mu_\mathrm{eq}$,
there exists a two-dimensional invariant slow manifold. The slow
coordinates $y$ and $\mu$ satisfy to first order the
equation}
\begin{equation}\label{eq:zie:system}
  \begin{aligned}
    \dot y=&\lambda_\mathrm{c}[y-y_\mathrm{eq}]+
    w_\mathrm{c}^\tran{}f_\mu [\mu-\mu_\mathrm{ref}]+a w_\mathrm{c}^\tran{}f_u\left[K_{\mathrm{st},y}[y-y_\mathrm{ref}]+K_{\mathrm{st},\mu}[\mu-\mu_\mathrm{ref}]\right]\\
    &+O(\|(y-y_\mathrm{eq},\mu-\mu_\mathrm{eq})\|^2)\mbox{,}\\
    \dot \mu=&K_{\mathrm{st},y}[y-y_\mathrm{ref}]+K_{\mathrm{st},\mu}[\mu-\mu_\mathrm{ref}]\mbox{,}
  \end{aligned}
\end{equation}
which has the Jacobian in the equilibrium
$(x_\mathrm{eq},\mu_\mathrm{eq})$
\begin{align*}
  A_\mathrm{ZIE}=
  \begin{pmatrix*}[c]
    \lambda_\mathrm{c}+a w_\mathrm{c}^\tran{}f_uK_{\mathrm{st},y}&
    w_\mathrm{c}^\tran{}f_\mu+a w_\mathrm{c}^\tran{}f_uK_{\mathrm{st},\mu} \\
    K_{\mathrm{st},y} & K_{\mathrm{st},\mu}
  \end{pmatrix*}\mbox{.}      
\end{align*}
Thus, criteria for stabilizing gains are that
\begin{align}\label{1d:combi:gaincrit1}
 K_{\mathrm{st},\mu}+a w_\mathrm{c}^\tran{}f_uK_{\mathrm{st},y}&<-\lambda_\mathrm{c}\mbox{,}&
\lambda_\mathrm{c}K_{\mathrm{st},\mu}-w_\mathrm{c}^\tran{}f_\mu K_{\mathrm{st},y}&>0\mbox{.}
\end{align}
Consequently, a sufficient non-degeneracy condition for the existence of stabilizing gains $(a,K_{\mathrm{st},y},K_{\mathrm{st},\mu})$ is that
\begin{align}
  \label{1d:combi:gaincrit2}
  w_\mathrm{c}^\tran{}f_\mu&\neq0\mbox{\qquad or\qquad}
  \left(\lambda_\mathrm{c}\neq0\mbox{\quad and\quad}w_\mathrm{c}^\tran{}f_u\neq0\right)
  \mbox{\qquad or\qquad}\lambda_\mathrm{c}<0\mbox{,}
\end{align}
which is violated only at events of codimension $2$ (if
$\lambda_\mathrm{c}<0$, $K_{\mathrm{st},y}=K_{\mathrm{st},\mu}=0$
stabilizes such that the condition $w_\mathrm{c}^\tran{}f_u\neq0$
is not necessary in that case).  If one of the first two cases of \eqref{1d:combi:gaincrit2} is
satisfied, we can adjust the input scaling $a$ such that the vectors
$(1,aw_\mathrm{c}^\tran{}f_u)$ and
$(\lambda_c,-w_\mathrm{c}^\tran{}f_\mu)$ are linearly
independent. Then the gains $(K_{\mathrm{st},y},K_{\mathrm{st},\mu})$
can be chosen from the quadrant defined by the two affine inequalities
in \eqref{1d:combi:gaincrit1}. The second condition on the gains in
\eqref{1d:combi:gaincrit1} is identical to the condition for control
through only the bifurcation parameter, \eqref{1d:SP}. Thus, it can be
approximated by a geometric condition, such that we obtain the
approximate criteria
\begin{align}\label{1d:combi:gaincrit3}
 K_{\mathrm{st},\mu}+a w_\mathrm{c}^\tran{}f_uK_{\mathrm{st},y}&<-\lambda_\mathrm{c}\mbox{,}&
\sigma\left[\partial_s\mu_\mathrm{eq}K_{\mathrm{st},\mu}+\partial_sy_\mathrm{eq}K_{\mathrm{st},y}\right]+O(\epsilon)&>0\mbox{.}
\end{align}
\begin{figure}[t]
\centering
\label{fig:flow_manipulation_js}
\includegraphics[width=\textwidth]{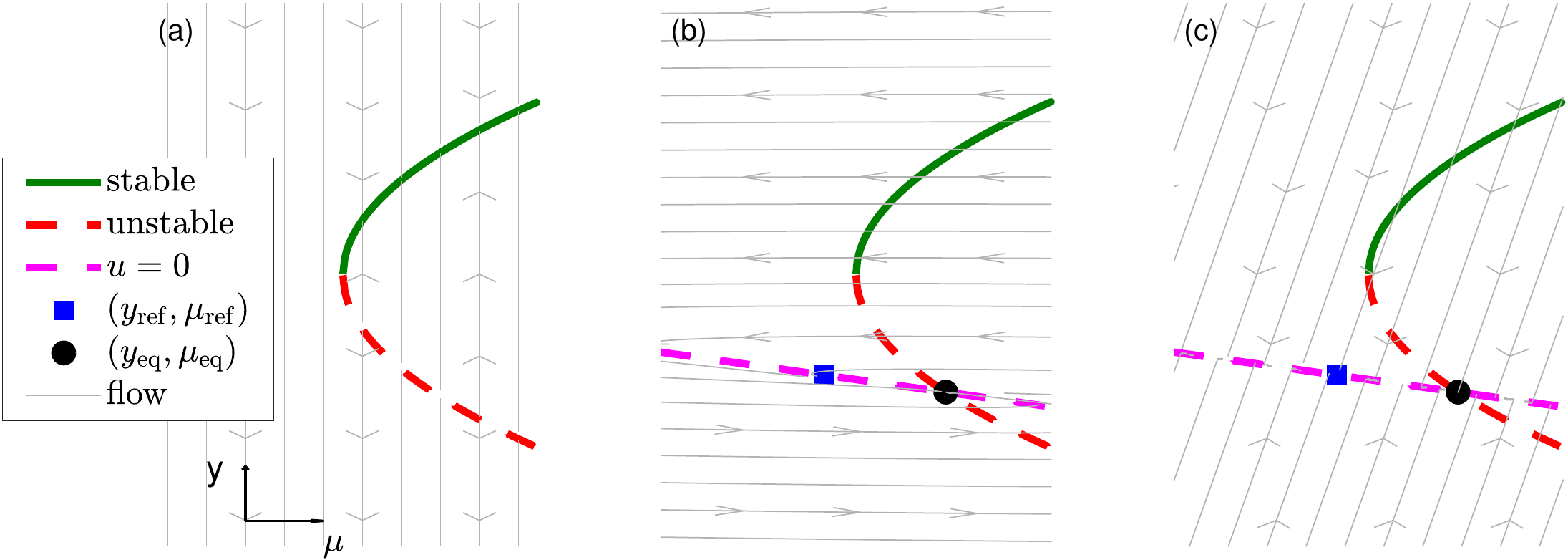}
\caption{\js{Illustration of the effect of different feedback control schemes on the flow. The stable (green) and unstable (red dashed) part of the equilibrium branch $(y_\mathrm{eq}(s),\mu_\mathrm{eq}(s))$ with the line $\{(y,\mu):u=K_{\mathrm{st},y}(y-y_\mathrm{ref})+K_{\mathrm{st},\mu}(\mu-\mu_\mathrm{ref})=0\}$ (dashed magenta) intersect in an equilibrium $(y_\mathrm{eq},\mu_\mathrm{eq})$ (black circle) for suitable gains $(K_{\mathrm{st},y},K_{\mathrm{st},\mu})$ and reference point $(y_\mathrm{ref},\mu_\mathrm{ref})$ (blue square).  Panel \textup{(}a\textup{):} uncontrolled system. Panel \textup{(}b\textup{):} control through bifurcation parameter in the case $|K_{\mathrm{st},\mu}|,|K_{\mathrm{st},y}|\gg1$ used by \cite{BartonSieber2013} (grey arrows are fast). Panel \textup{(}c\textup{):} zero-in-equilibrium control in the case $|a|\gg1$ (grey arrows are fast). The horizontal axis is $\mu$ and the vertical axis is $y$ as is convention for bifurcation diagrams.} %\jsq{removed all discussion, which should be in main text}%Especially in the case where the solution curve has a very steep fold( in other words if $|w_\mathrm{c}^\tran{}f_\mu|<<1$), then the line $u=0$ has a very large gain $K_{\mathrm{st},y}$. This may lead to large fluctuations in the y direction and, thus,  failure to stabilize the system. On the other hand, the presence of the control gain $a$ in the  zero-in-equilibrium feedback allows to stabilize the system with a gain $K_{\mathrm{st},y}$ of arbitrary small modulus.
}
\end{figure}

\js{\Cref{fig:flow_manipulation_js} illustrates how \js{feedback control through the bifurcation parameter and} zero-in-equilibrium affect the flow to stabilize unstable equilibria. \Cref{fig:flow_manipulation_js}(a) shows an equilibrium branch $(y_\mathrm{eq}(s),\mu_\mathrm{eq}(s))$ with saddle-node bifurcation without control (System \eqref{eq:zie:system} with $K_{\mathrm{st},y}=K_{\mathrm{st},\mu}=0$). The flow is indicated by grey arrows.}

\js{  \Cref{fig:flow_manipulation_js}(b) shows the effect of control through bifurcation parameter, which is a special case of zero-in-equilibrium feedback for $a=0$. The sketch shows the case of large gains ($|K_{\mathrm{st},y}|,|K_{\mathrm{st},\mu}|\gg1$), which approximates the feedback rule \eqref{eq:muBS}, $\mu(t)=\mu_\mathrm{ref}+K_\mathrm{st}(y-y_\mathrm{ref})$. The large gains cause the controlled flow to be slow-fast. In \cref{fig:flow_manipulation_js}(b,c) the grew arrows indicate the direction of the fast flow. The slow flow (not indicated in \cref{fig:flow_manipulation_js}(b,c)) then follows the line $\{u=0\}$ (dashed, magenta) toward the equilibrium $(y_\mathrm{eq},\mu_\mathrm{eq})$ (black circle). For the large-gain regime the fast flow is nearly horizontal (exclusively changing $\mu$), causing large corrections in $\mu$. For smaller gains $K_{\mathrm{st},y},K_{\mathrm{st},\mu}$ trajectories of the controlled system may spiral, crossing the line $\{u=0\}$ (dashed magenta), vertically, but eventually (for suitable gains) converging to the equilibrium $(y_\mathrm{eq},\mu_\mathrm{eq})$ (black circle) which lies in the intersection of the line $u=0$ and the equilibrium curve.}

\js{  \Cref{fig:flow_manipulation_js}(c) shows the effect of zero-in-equilibrium control ($a\neq0$), emphasizing its effect by choosing $|a|\gg1$. Again, the large parameter makes the controlled flow slow-fast. However, for zero-in-equilibrum control the fast flow is nearly parallel to the lines $\{\mu=\mathrm{const}\}$. Thus, the control gain $a$ enables control without large deviations in the parameter $\mu$. This is beneficial for pedestrian flow control in a physical experiment, where the parameter is the location of the obstacle where large parameter variations correspond to large-amplitude motions of the obstacle with real-time requirements. \Cref{fig:ped_line} in \cref{sec:res_zie} also shows that even for the pedestrian flow simulation the required size of gains for control through the bifurcation parameter pushes the simulation out of the regime where one may assume that there is only one slow dimension (see \cref{ass:1d}). % , governs the flow manipulation leading to a more robust control scheme which leads to convergence even in the presence of noise.
} 

Criterion \eqref{1d:combi:gaincrit1} \js{expresses this  advantage of
zero-in-equilibrium control %control through the extended washout filters 
over 
control purely through the bifurcation parameter $\mu$ (corresponding
to $a=0$) quantitatively.}
% , beyond the avoidance of loss of stabilizability at codimenion-$1$ events.
If the coefficient
$w_\mathrm{c}^\tran{}f_\mu$ is relatively small compared to
$\lambda_\mathrm{c}>0$ (so the system does not react quickly to
changes in $\mu$ along an unstable branch), the corrections by
feedback control through $\mu$ (in $\dot \mu=u$) have to be large when
$a=0$, because $K_{\mathrm{st},\mu}<-\lambda_\mathrm{c}$ is required when $a=0$, which
implies $\lambda_\mathrm{c}K_{\mathrm{st},\mu}<-\lambda_\mathrm{c}^2$,
such that
$|K_{\mathrm{st},y}|>\lambda_\mathrm{c}^2/|w_\mathrm{c}^\tran{}f_\mu|$
is required, which can be large, resulting in large right-hand sides
for $\dot\mu=u$ in the presence of small disturbances. \js{Geometrically this means that the line $\{u=0\}$ in \cref{fig:flow_manipulation_js}(b) would be almost horizontal, resulting in trajectories of the controlled flow that are also almost horizontal (grey lines in \cref{fig:flow_manipulation_js}(b)) leading to large corrections  in $\mu$. % and prone to noise and small disturbances, leading possibly to failure to stabilize the corresponding equilibrium
 On the other hand,} the additional
input $au$ in the right-hand side of $f$ with a (possibly large) scaling $a$
of suitable sign permits us to satisfy the first criterion in
\eqref{1d:combi:gaincrit1} with a gain $K_{\mathrm{st},y}$ of
arbitrary small modulus and $K_{\mathrm{st},\mu}>0$, such that the
second criterion is also satisfied (thus, most control is exerted
through the input $u$ in the right-hand side $f$). \js{\Cref{fig:flow_manipulation_js}(c) shows that even for a nearly horizontal line $\{u=0\}$ (required if $|w_\mathrm{c}^\tran{}f_\mu|\ll1$) the controlled flow does not show large excursions in the parameter $\mu$, converging well to its equilibrium. }

The particle-flow model for the pedestrian evacuation scenario
introduced in \Cref{sec:scenario} has this feature: controlling the
flow by shifting the obstacle in $x$-direction requires large control
action to compensate small disturbances. The additional input applies
to each pedestrian a biasing force, with a strong direct effect on the
output flux measures, such that the relative weighting between the
control inputs is $a=50$ (see \Cref{sec:gains_choice}).

\section{Pedestrian Evacuation Scenario}
\label{sec:scenario}
In this section we describe a prototype multi-particle model for an evacuation scenario. We refer to it as a microscopic model, as we set the rules for the dynamics at the level of individual pedestrians. As mentioned in \Cref{sec:intro:pedestrian},  the system exhibits tipping, a sudden change from one stable state to another stable state, and hysteresis, which we will analyze with non-invasive feedback control. We will treat the microscopic model like a physical experiment in the sense that we assume similar limitations. In particular,  the microscopic state of the system is not set 'at will' and  precise derivative information is not available.
\begin{figure} [ht]
\centering
\includegraphics[width=0.7\textwidth]{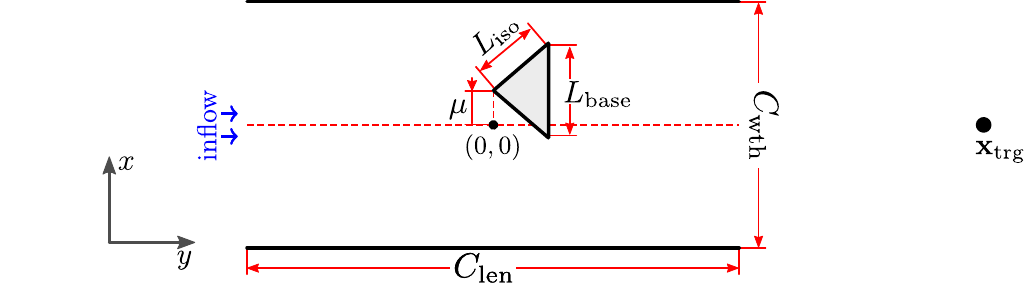}
\caption{Geometry of corridor, obstacle and inflow area of pedestrians. See \Cref{tab:model:parameters} for values of the parameters for corridor and obstacle geometry.}
\label{fig:corridor}
\end{figure}

\subsection{General Set-up}
It is assumed that pedestrians want to evacuate a building, passing through a corridor with an obstacle, as shown in \cref{fig:corridor}. To exit the building via the corridor, pedestrians have to choose which route they want to follow to maneuver around the obstacle. This choice is influenced by the shortest way to the exit and by the walking behavior of nearby pedestrians. The route choice behavior is investigated by changing the position of the obstacle and, thus, the preference for each route.

\subsubsection*{Social force model with lemming effect}
\begin{figure} [!htb]
\centering
\subfloat[Graphs of  repulsive force $F^s_\mathrm{rep}(r)$ in \eqref{eq:pedped}, for $s\in\{\mathrm{ped},\mathrm{obj}\}$. The function $F^s_\mathrm{rep}(r)$ is zero for $r \geq \sigma^s$.]{\label{fig:repforce}\includegraphics[width=0.47\textwidth]{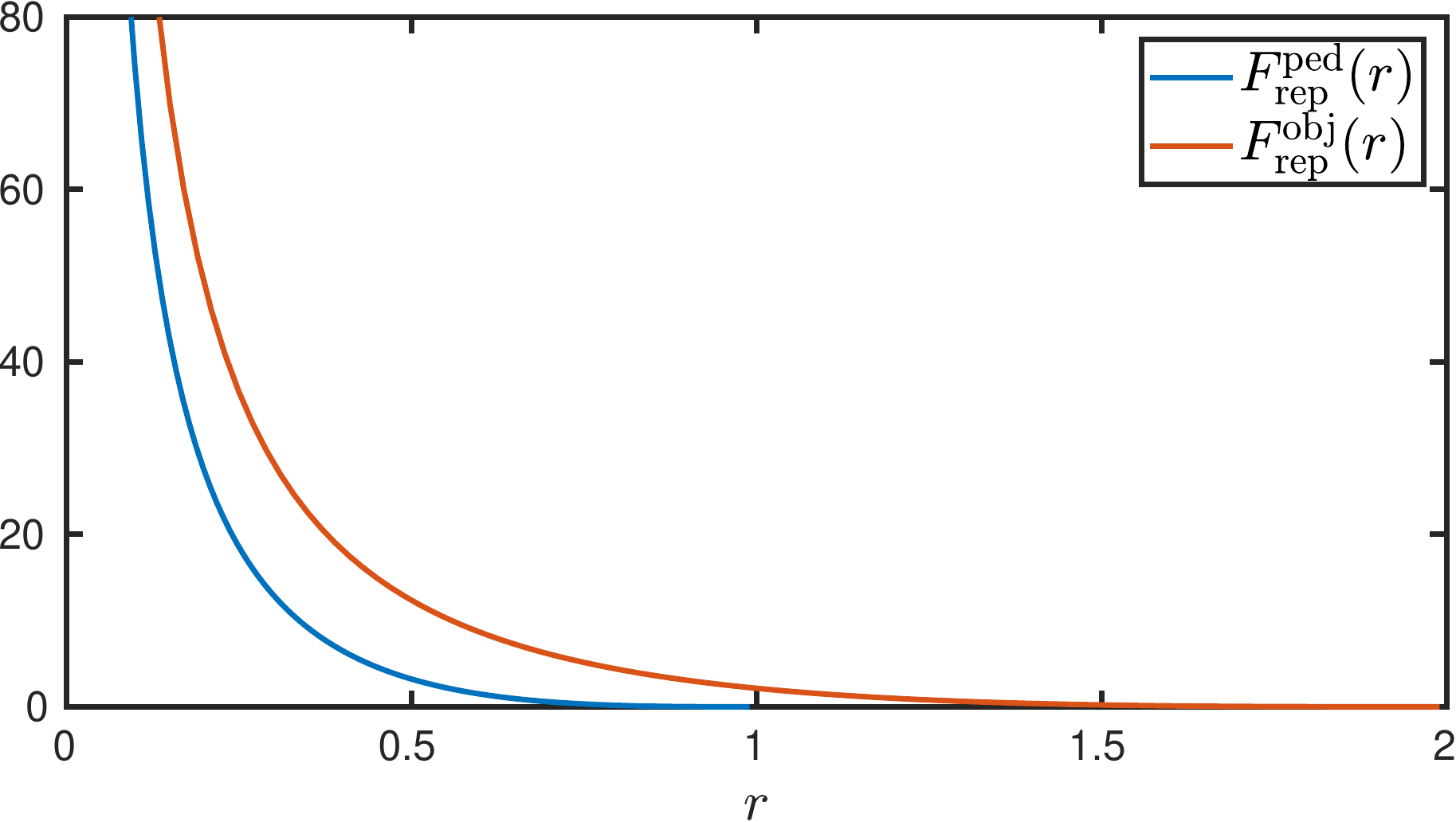}}
\hfill
\subfloat[\footnotesize Alignment weighting $\kappa(r,\theta)$ in \eqref{eq:kappa}, shown as a graph of the complex argument $r\exp(\mathrm{i}\theta)$.]{\label{fig:kappa}\includegraphics[width=0.47\textwidth]{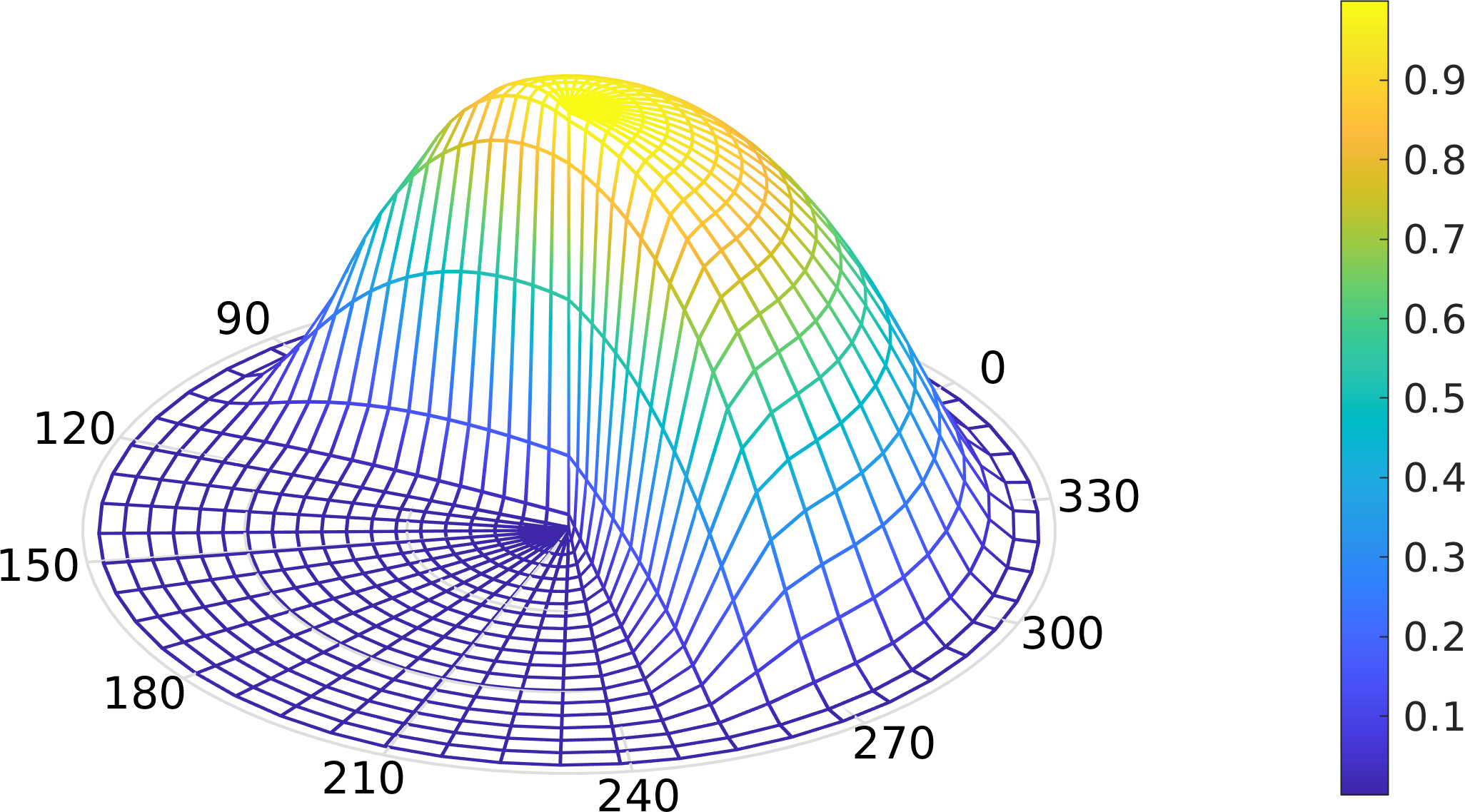}}
\caption{Graphs of interaction forces for repulsion \textup{(}\Cref{fig:repforce}\textup{)} and alignment \textup{(}\Cref{fig:kappa}\textup{)}. In \Cref{fig:kappa} the radial component is the distance of pedestrian $j$ from $i$, the angular component is the difference of their  angular velocity. See \Cref{tab:model:parameters} for parameter values.}
\label{fig:dddd}
\end{figure}

We consider a scenario with $N$ pedestrians, where $N$ is large.
Helbing and Molnar proposed a model in \cite{HelbingSocial} where each
pedestrian $i$ is described by a particle of zero extent at
position $\mathbf{x}_i(t)$ and moving with velocity
$\Dot{\mathbf{x}}_i(t)$ in the plane. Their motion is the response to
forces acting on it. The so-called \emph{social force model} assumes
that the main forces that determine the motion of pedestrian $i$ are
their tendencies to do the following.
\begin{itemize}
\item Pedestrian $i$ moves  towards a target point $\mathbf{x}_\mathrm{trg}\in\R^2$ (see \Cref{fig:corridor}) aiming for a desired speed $v_\mathrm{trg}$. This results in a target attraction force $\mathbf{F}_{\mathrm{trg},i}$ acting on pedestrian $i$.
\item Pedestrian $i$ avoids close encounters with pedestrian $j$ (for
  all $j\neq i$, $j\leq N$), resulting in a repulsive force $\mathbf{F}_{\mathrm{rep},ij}^\mathrm{ped}$.
\item  Pedestrian $i$ avoids  collision with each object $j$, which can be an obstacle or wall, resulting in a repulsive force $\mathbf{F}_{\mathrm{rep},ij}^\mathrm{obj}$.
\end{itemize} 
The target attraction force is
\begin{align}
\label{eq:dir_force}
\mathbf{F}_{\mathrm{trg},i}&= \frac{1}{\tau}(v_\mathrm{trg}\mathbf{e}_{\mathrm{trg},i} - \mathbf{\dot x}(t))\mbox{,\quad where\quad}
\mathbf{e}_{\mathrm{trg},i}=\frac{\mathbf{x}_\mathrm{trg}-\mathbf{x}_i}{\| \mathbf{x}_\mathrm{trg}-\mathbf{x}_i\|}
\end{align}
is the direction vector toward the target point
$\mathbf{x}_\mathrm{trg}$ and $\tau $ is the reaction time. In our
case all pedestrians have the same target point and reaction time (see
\Cref{tab:model:parameters}). The force $\mathbf{F}_{\mathrm{trg},i}$ such
that, in the absence of other pedestrians or obstacles, pedestrian $i$
adjusts their velocity with rate $\tau$ such that they move with speed
$v_\mathrm{trg}$ toward $\mathrm{x}_\mathrm{trg}$. The repulsive
forces from other pedestrians and from objects are modeled as
monotonic decreasing functions of the distance to other pedestrians
and obstacles respectively. Following \cite{molnar2001control}, we make the following choice for
pedestrian-pedestrian (superscript $s=\mathrm{ped}$) and
pedestrian-obstacle/wall (superscript $s=\mathrm{obj}$)
interactions forces:
\begin{align}\allowdisplaybreaks
\label{eq:pedped}
\mathbf{F}_{\mathrm{rep},ij}^s&=F_\mathrm{rep}^s(\|\mathbf{r}_{ij}\|)\frac{\mathbf{r}_{ij}}{\|\mathbf{r}_{ij}\|}\mbox{,\ with\ }
F_\mathrm{rep}^s(r)=
\begin{cases}
  -V_\mathrm{rep}^s\left[\tan g^s(r)-g^s(r) \right]&\mbox{if $r<\sigma^s$,}\\
  0&\mbox{if $r\geq\sigma^s$,}
\end{cases}\mbox{ and}\\
\nonumber
g^s(r)&=\frac{\pi}{2}\left(\frac{r}{\sigma^s}-1\right)\mbox{\quad for $s=\mathrm{obj}$ or $\mathrm{ped}$.}
\end{align}
Here $\mathbf{r}_{ij}=\mathbf{x}_j-\mathbf{x}_i$ is the vector between
pedestrian $i$ and pedestrian (or obstacle) $j$. For the case of an
obstacle (wall or triangle), this vector is defined as pointing toward
the point of the obstacle $j$ that is closest to the pedestrian $i$. The parameters
$V_\mathrm{rep}^\mathrm{ped}$ and $V_\mathrm{rep}^\mathrm{obj}$
control the repulsion strength between pedestrians, or between
pedestrians and obstacles, respectively (chosen uniform for all
pedestrians and obstacles here). The repulsion force is radially symmetric and has finite range, in contrast to \cite{HelbingSocial}. \Cref{fig:repforce} shows its dependence on the distance $r=\|\mathbf{r}_{ij}\|$ between particles or objects $i$ and $j$. So, the equation of motion for
pedestrian $i\leq N$ according to the social force model is given by
\begin{align*}
\mathbf{\ddot x}_i= \mathbf{F}_{\mathrm{trg},i} + \sum_{j=1}^N\mathbf{F}_{\mathrm{rep},ij}^\mathrm{ped}
 + \sum_{\mathrm{objects\ }k} \mathbf{F}_{\mathrm{rep},ik}^\mathrm{obj}\mbox{\quad (social force model).}
\end{align*}

Starke \emph{et al}.~\cite{starkeieeeqingdao14} hypothesize
the presence of % adapted the model by taking into account
another force, namely a tendency to follow others. For example, this effect may be %is especially
present in an emergency situation when there is no good knowledge of the
geometry of the building. In such a situation pedestrian $i$ has a preference to move in the same direction as other pedestrian around
him. This psychological factor was called the \emph{lemming effect} in \cite{starkeieeeqingdao14}
and was modeled by changing the directional vector in $\mathbf{F}_{\mathrm{trg},i}$
to a linear combination of the direction $\mathbf{e}_{\mathrm{trg},i}$
towards the target point $\mathbf{x}_\mathrm{trg}$ and a weighted mean velocity
$\langle\textbf{v}\rangle_i$ of the velocity vectors
$\textbf{v}_j=\dot{\textbf{x}}_j$ of pedestrians $j$ in a
neighborhood of pedestrian $i$:
\begin{align}\label{model:eal}
\mathbf{e}_{\mathrm{al},i}(t)&=\cfrac{(1-p_\mathrm{al}) \mathbf{e}_{\mathrm{trg},i} +
  p_\mathrm{al}\langle\mathbf{v}\rangle_i}{\|(1-p_\mathrm{al}) \mathbf{e}_{\mathrm{trg},i} +
  p_\mathrm{al}\langle\mathbf{v}\rangle_i\|}\mbox{,\quad where}&
\langle\mathbf{v}\rangle_i&=\cfrac{\sum_{j\neq i} \kappa(\|\mathbf{r}_{ij}\|,\angle\,\Dot{\mathbf{r}}_{ij}) \mathbf{v}_j(t)}{\sum_{j\neq i} \kappa(\|\mathbf{r}_{ij}\|,\angle\,\Dot{\mathbf{r}}_{ij})}\mbox{,}
\end{align}
and $p_\mathrm{al} \in [0,1]$  denotes the \emph{lemming parameter} which controls the influence other pedestrians have over the target direction. The weight function $\kappa(r,\theta)$ depends on the
distance $r=\| \mathbf{r}_{ij} \|$ of pedestrians $i$ and $j$ and the
angle
$\theta=\angle\,\Dot{\mathbf{r}}_{ij}=\angle\Dot{\mathbf{x}}_j-\angle\,\Dot{\mathbf{x}}_j$
between their walking directions (their velocity vectors).
The quantity $\langle\textbf{v}\rangle_i$ is called the weighted mean velocity vector. It depends on weights determined by the real-valued weighting function
\begin{align}
\label{eq:kappa}
\kappa(r,\theta)&=
\begin{cases}
  \cfrac{\gamma}{1+\exp(-\alpha \cos (\beta \theta))}\,\exp\left(\cfrac{\sigma_\mathrm{al}^2}{r^2-\sigma_\mathrm{al}^2}\right)&\mbox{if $r\leq \sigma_\mathrm{al}$,}\\
  0&\mbox{if $r>\sigma_\mathrm{al}$}
\end{cases}\mbox{\quad for $r\geq0$, $\theta \in[-\pi,\pi]$.}
\end{align}
   The weight
function $\kappa$ has a finite support radius $\sigma_\mathrm{al}$ for
$r$. The parameter $\gamma$ is a scaling factor, and parameters
$\alpha$ and $\beta$ are chosen such that each pedestrian is
influenced by pedestrians nearby, walking in the same direction. More
precisely, each pedestrian is influenced by others walking in
approximately the same directions ($\kappa$ is noticeably positive for angles
$\theta \in [-100^\circ, 100^\circ]$ and for
$r<\sigma_\mathrm{al}$, see \Cref{fig:kappa}). Assuming that pedestrian $i$ is
placed in the center, the graph in \Cref{fig:kappa} indicates how much
another pedestrian $j$ inside the finite support radius influences the
mean $\langle \mathbf{v}\rangle_i$ depending on their relative walking
direction.

Consequently, the target attraction force is modified to take into account the tendency for alignment, such that we modify our social force model. The equation of motion for pedestrian $i$ is given by 
\begin{align}\label{model:xddot}
\mathbf{\ddot x}_i&= \mathbf{F}_{\mathrm{al},i} + \sum_{j=1}^N\mathbf{F}_{\mathrm{rep},ij}^\mathrm{ped}
+ \sum_{\mathrm{objects\ }k} \mathbf{F}_{\mathrm{rep},ik}^\mathrm{obj}\mbox{\quad (social force model with alignment), where}\\
\nonumber
\mathbf{F}_{\mathrm{al},i}&= \frac{1}{\tau}(v_\mathrm{trg}\mathbf{e}_{\mathrm{al},i} - \mathbf{\dot x}(t))\mbox{,}
\end{align}
and $\mathbf{e}_{\mathrm{al},i}$ is defined in \eqref{model:eal}, and
$\mathbf{F}_{\mathrm{rep},ij}^\mathrm{ped}$ and
$\mathbf{F}_{\mathrm{rep},ik}^\mathrm{obj}$ are defined in \eqref{eq:pedped}.
\subsection{Bistable behavior and hysteresis}
In the following, we consider model \eqref{model:xddot}  with $N=100$ pedestrians in the corridor. The boundary conditions are such that for every pedestrian exiting the corridor at $y=C_\mathrm{len}/2$, 
a new pedestrian enters at $y=-C_\mathrm{len}/2$ with initial velocity $(v_\mathrm{trg},0)$ and with a vertical position uniformly randomly distributed around the center line of the corridor within the interval $[-0.5,0.5]$. The values of the remaining parameters can be found in Table \ref{tab:model:parameters}. The number of pedestrians is chosen so that the crowd in the corridor is of medium density, too many pedestrians would overcrowd the corridor while too little would result in an interrupted flow. The system parameter which we vary to perform the bifurcation analysis is the position of the tip of the triangular obstacle $\mu$ as shown in \Cref{fig:corridor}.

\paragraph{Details of simulation protocol for \Cref{fig:intro}(b)} Starting with  $\mu=-1.2$\,m, the system was numerically integrated (see \Cref{app:sec:parameters} for details of integration). The parameter $\mu$ was increased in steps of $0.1$ meters every 300 seconds, so we were slowly changing the position of the tip of the triangle performing a quasi-stationary \emph{up-sweep}. When $\mu=1.2$ m, we started decreasing it in steps of $0.1$ m until the triangle was at its original position ($\mu=-1.2 $ m) performing a \emph{down-sweep}.

Reproducing the results of \cite{starkeieeeqingdao14}, we 
confirmed that the system exhibits bistability and hysteresis, as
shown in \Cref{fig:intro}(b). The macroscopic variable used for the
$y$-axis in \Cref{fig:intro}(b) is the difference of fluxes
$\Delta\Phi=\Phi_+-\Phi_-$ at each side of the obstacle.  In
\cref{fig:intro}(b), the fluxes $\Phi_+$ and $\Phi_-$ are measured as
the number of pedestrians passing the end of the obstacle. More
precisely,
\begin{align}
  \label{eq:DPhi}
  \Delta\Phi(t)&=\frac{1}{\tau_{\max}}\int_{t-\tau_{\max}}^t \bigl( \Phi_+(s)-\Phi_-(s) \bigr)\,\mathrm{d} s\mbox{,\quad where}\\
  \nonumber
  \Phi_\pm(t)&=\sum_{i=1}^NE_\Phi^\pm(\mathbf{x}_i(t))\mbox{,\ and\ }
  E_\Phi^\pm(x,y)=
  \begin{cases}
    1&\mbox{if $\pm(x-\mu)>0$ and $|y-y_{\mathrm{c},\Phi}|\leq y_{\mathrm{len},\Phi}$,}\\
    0&\mbox{otherwise.}
  \end{cases}
\end{align}
Thus, $\Phi_\pm(t)$ is the number of pedestrians with position
$\mathbf{x}_i(t)=(x_i(t),y_i(t))$ in a spatial box given by
$|y_i-y_{\mathrm{c},\Phi}|\leq y_{\mathrm{len},\Phi}$ and $x_i>\mu$
for $\Phi_+$ or $x_i<\mu$ for $\Phi_-$  (see also \cref{fig:pedweight}). The concrete parameters for
our spatial box for counting are
$y_{\mathrm{c},\Phi}=\sqrt{\smash[b]{L_\mathrm{iso}^2-L_\mathrm{base}^2/4}}+0.5\,\mathrm{m}=2.5$\,m,
$y_{\mathrm{len},\Phi}=0.5$\,m. The
averaging interval length $\tau_{\max}$ for $\Delta\Phi(t)$ is $10$
seconds. In \cref{fig:intro}(b), the value of
$\Delta\Phi=\Phi_+-\Phi_-$ at the end of the $300$\,s interval, just
before we change the position of the obstacle $\mu$, is plotted for
each value of $\mu$. %}}
The two overlapping stable branches of steady flow states in
\Cref{fig:intro}(b) are the result of this up-sweep and down-sweep of
the parameter $\mu$. At one value of $\mu$ for each sweep, the steady
flow state jumps from one branch to the other. This hysteresis suggests the
existence of an unstable branch of steady states that connects the
stable ones at two saddle-node bifurcations.

\subsection{Output and input for measurement and control}
\label{sec:outin}

Due to the randomness of the $x$-coordinate at entry to the corridor, we can expect deterministic results in our system only in the limit number of pedestrians $N\to\infty$ with suitably scaled corridor and obstacle parameters $\mu,L_\mathrm{base},L_\mathrm{iso},C_\mathrm{len},C_\mathrm{wth}\sim\sqrt{N}$.  The fluctuations due to finite $N$ make the averaging over time for the fluxes $\Phi_\pm$ necessary.
The flux measure $\Phi_\pm$ has further disadvantages, especially if one plans to make a feedback control input depend on it  to perform bifurcation analysis. The measure can only take a finite and small number of integer values. In addition, the results depend on the time window size for the average, where one has a tradeoff. A small time window results in large fluctuations, while a large window averages out important system characteristics and may introduce a delay into the control feedback loop.
\subsubsection{Instantaneous space-averaged flux measure}
We introduce a space-averaged flux measure $\phi$, which is
instantaneous. This means that the measure $\phi(t)$ depends only on the vector
$(\mathbf{x}(t),\Dot{\mathbf{x}}(t))$, not a history
$\tau\mapsto(\mathbf{x}(t+\tau),\Dot{\mathbf{x}}(t+\tau))$ for $\tau$
in some interval $[-\tau_\mathrm{max},0]$ (in contrast to
$\Phi_\pm(t)$, where $\tau_\mathrm{max}=10$\,s).

The space average has a weight kernel $w_\mathrm{pos}(x,y)$ such that every pedestrian counts with a weight according to their position. The contribution of each pedestrian to the flux is the product of this value with its velocity component $v^x_i(t)$ along the corridor. 
\begin{figure}[t]
\centering
\subfloat[Control signal $u$ acts as bias force on pedestrians in red box in front of obstacle, where $b(x,y)\neq 0$ in \eqref{eq:box}.]{\label{fig:box}\raisebox{4.7ex}{\includegraphics[width=0.38\textwidth]{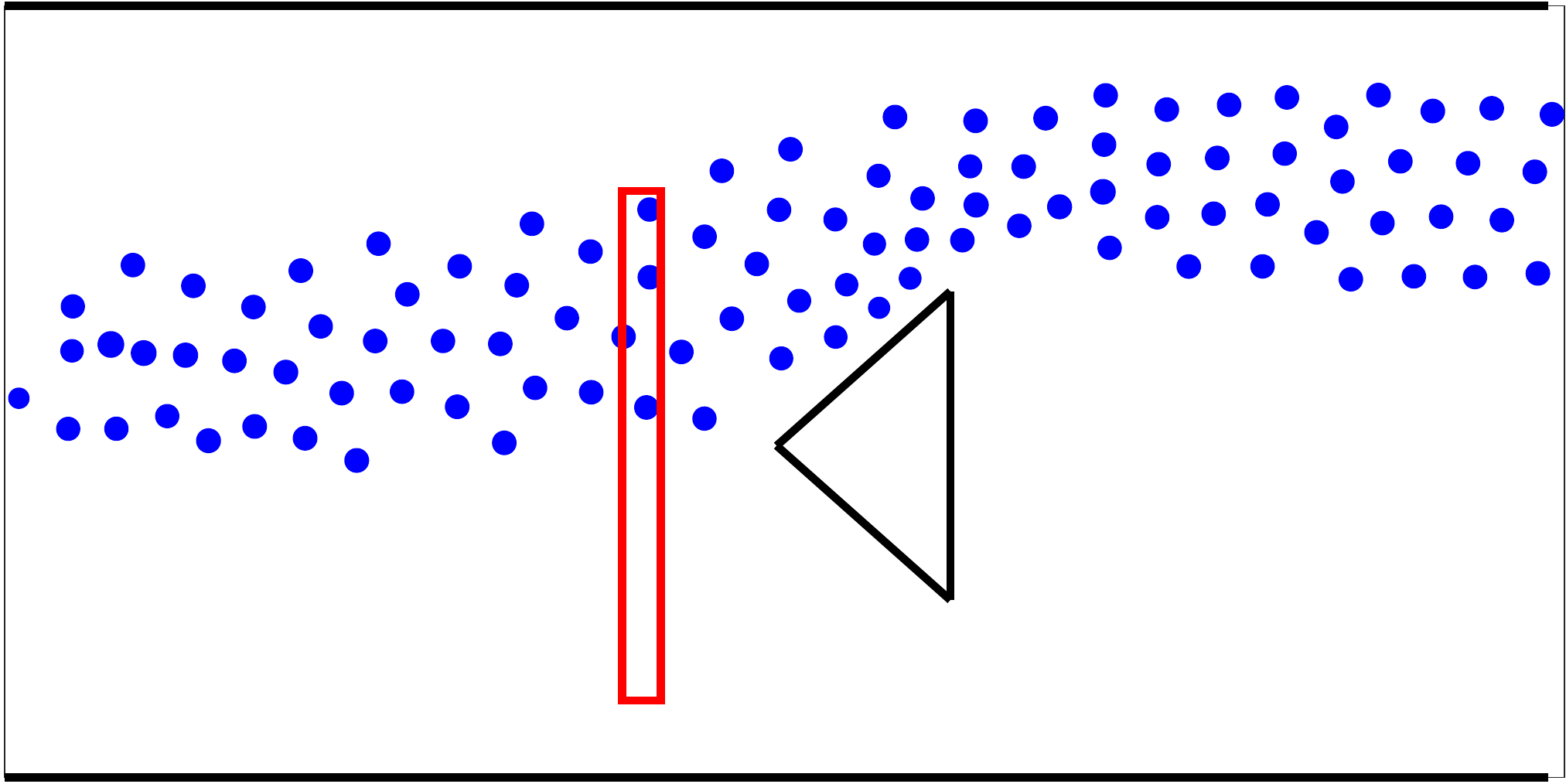}}}
\hspace*{1em}
\subfloat[Color: weight $w_\mathrm{pos}$ for pedestrian positions in space averaged flux measure $\phi$; Dashed and dotted boxes: area where $E^\pm_\Phi\neq 0$ in \eqref{eq:DPhi}.]{\label{fig:pedweight}\includegraphics[width=0.49\textwidth]{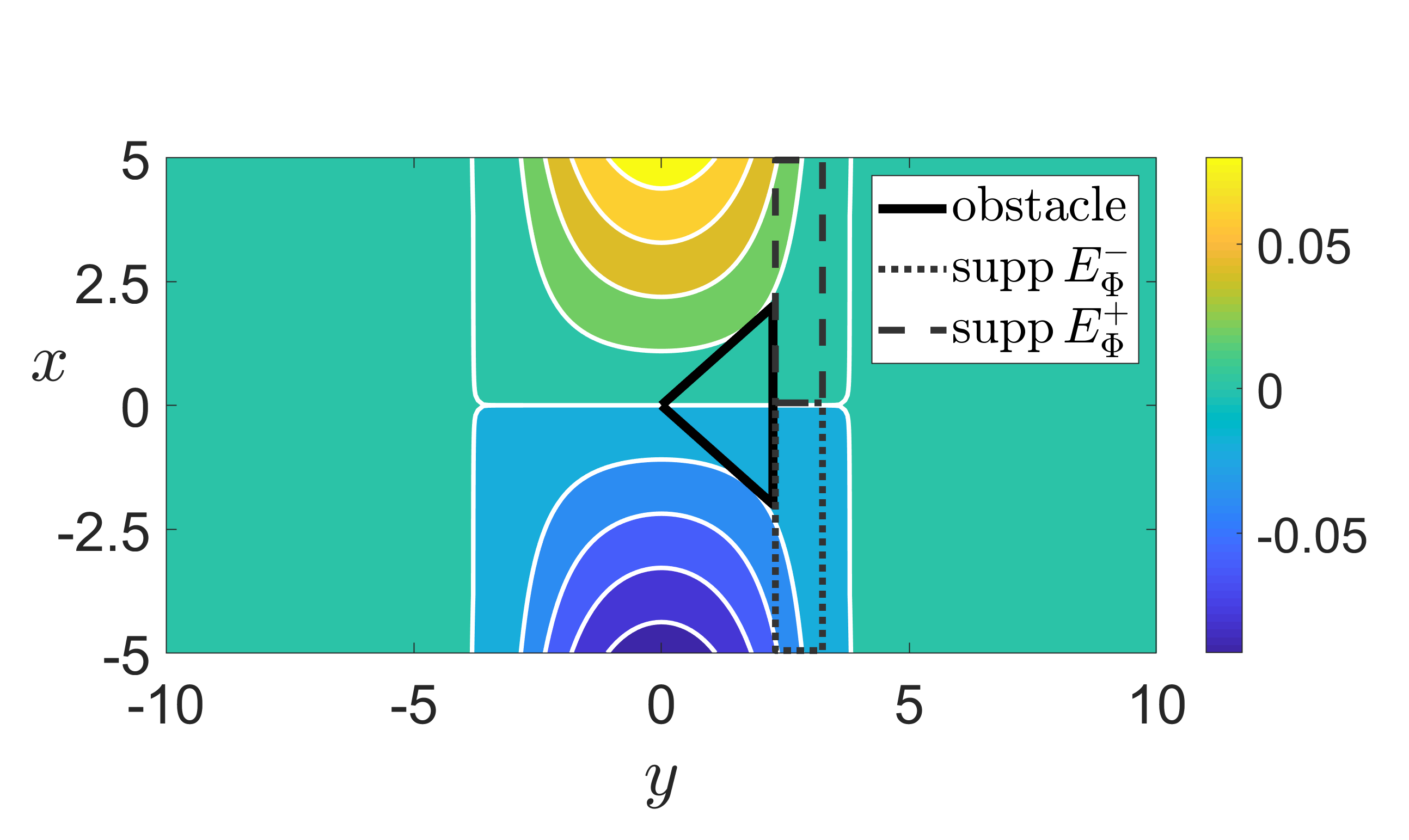}}
\unitlength 0.8cm
%\put(-5.6,0.2){$y$}
%\put(-9.2,3.6){$x$}
\caption{Area for input \textup{(}\Cref{fig:box}\textup{)}, and weight
  function $w_\mathrm{pos}(x,y)$ for instantaneous flux measure, given
  in \eqref{w}. At $y=0$, $w_\mathrm{pos}(\cdot,0)$ decreases linearly
  in $x$ (maximal at upper wall, here $x=5$, minimal at lower,
  $x=-5$. \Cref{fig:pedweight} also shows the boxes
  $E_\Phi^\pm(x,y)$ used for counting in flux measure $\Phi_\pm$ in
  \eqref{eq:DPhi}.}
\label{fig:inout}
\end{figure}
Denote the state space vector for pedestrian $i$ by $\mathbf{X}_i(t)=(\mathbf{x}_i(t),\Dot{\mathbf{x}}_i(t))=(x_i(t) ,y_i(t),v^{x}_i(t),v^y_i(t))$, and the overall state space vector by $\mathbf{X}(t)=(\mathbf{X}_1(t),\ldots,\mathbf{X}_N(t))$ at time $t$.
Then the flux $\phi$ is defined as
\begin{align}\allowdisplaybreaks
%\begin{equation}
\phi(\mathbf{X}(t))&=\sum_{i=1}^N w_\mathrm{pos}(x_i(t),y_i(t))v_i^{x}(t)\mbox{, with}\label{flux}\\
%\end{equation}
w_\mathrm{pos}(x,y)&= E(x,r_+(x,y))-E(x,r_-(x,y)),\label{w}\
E(x,r)=
\begin{cases}
  \eta |x| \displaystyle\int\limits_{\frac{d^2}{d^2-r^2}}^{\infty}  \frac{\mathrm{e}^{-t}}{t}\,\mathrm{d} t&\mbox{if $|r| < d$}\\
  0&\mbox{if $|r| \geq d$,}
\end{cases}\\
r_\pm^2(x,y)&=(x\pm x_{\mathrm{c},\phi})^2+(y-y_{\mathrm{c},\phi})^2\mbox{.}\nonumber%\\ 
%r_2^2&=(x_i-c_1)^2+(y_i-c_2)^2, \nonumber
\end{align}
The parameter $d=4$\,m is the length scale of the weight function
$w_\mathrm{pos}$ and $\eta=(1/12)\,\mathrm{s}/\mathrm{m}^3$ is a scaling
factor which also makes the flux dimensionless. The weight function
was chosen to have two bell shaped humps (one positive, one negative, see \cref{fig:pedweight})
with extrema at the points $(\pm x_{\mathrm{c},\phi},y_{\mathrm{c},\phi})=(\pm C_\mathrm{wth}/2,0)$ where pedestrian motion should be weighted
highest: The $y$ coordinate of the extrema, $y_{\mathrm{c},\phi}$ is the horizontal ($y$)
position of the tip of the triangular obstacle ($y_{\mathrm{c},\phi}=0$), while the weight
increases linearly in $x$ along the line $y=0$.

\begin{figure}[t]
\centering
\subfloat[Bistability in the pedestrian model, shown using flux measure $\phi$ in the $(\mu,\phi)$-plane.]{\label{fig:fluxm_sweep}\includegraphics[height = 0.25\textheight]{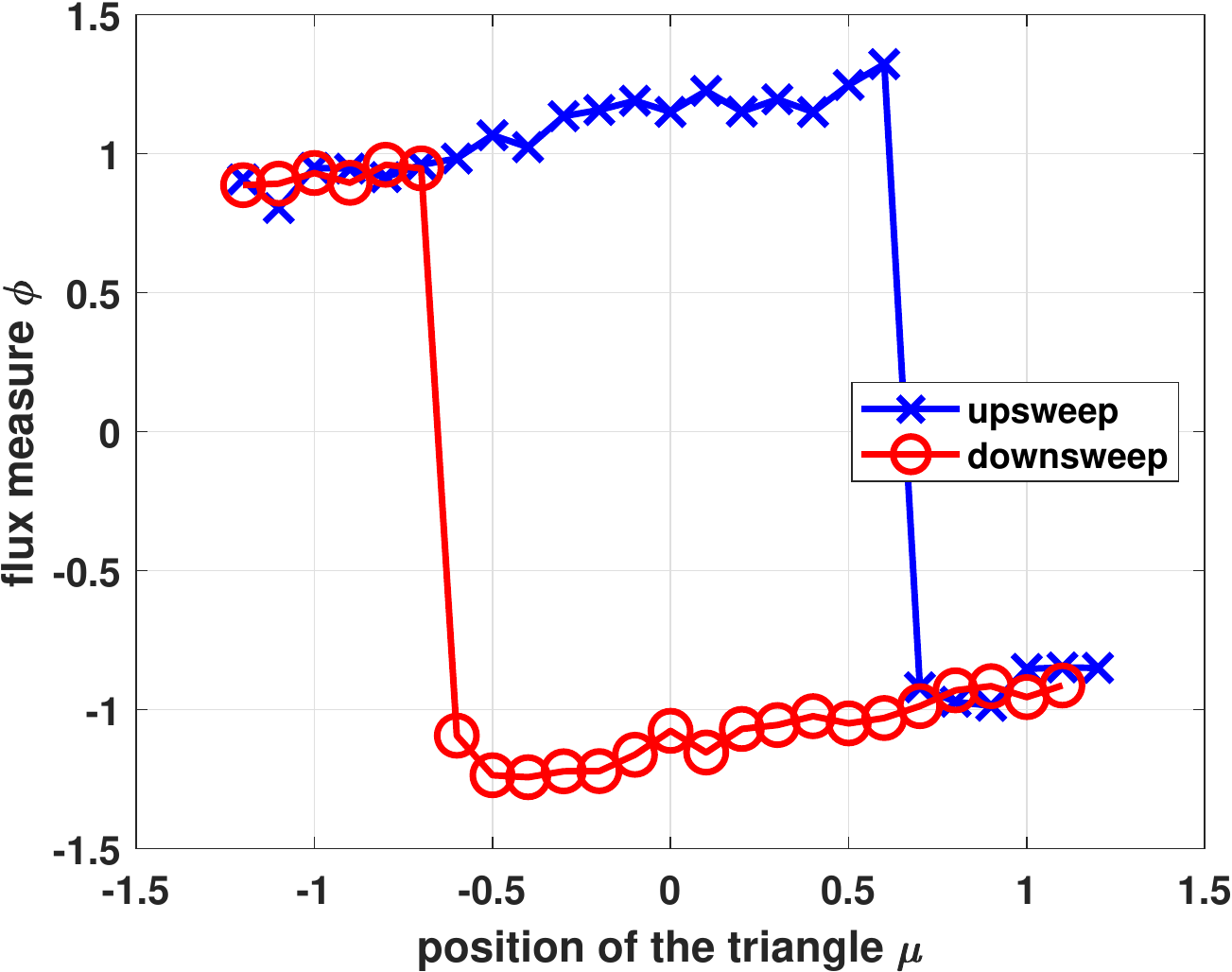}}\hspace*{1em}
\subfloat[Time series of the flux measure $\phi$ for  position of the obstacle $\mu=-1.2$ m.]{\label{fig:fovert}\includegraphics[height = 0.25\textheight]{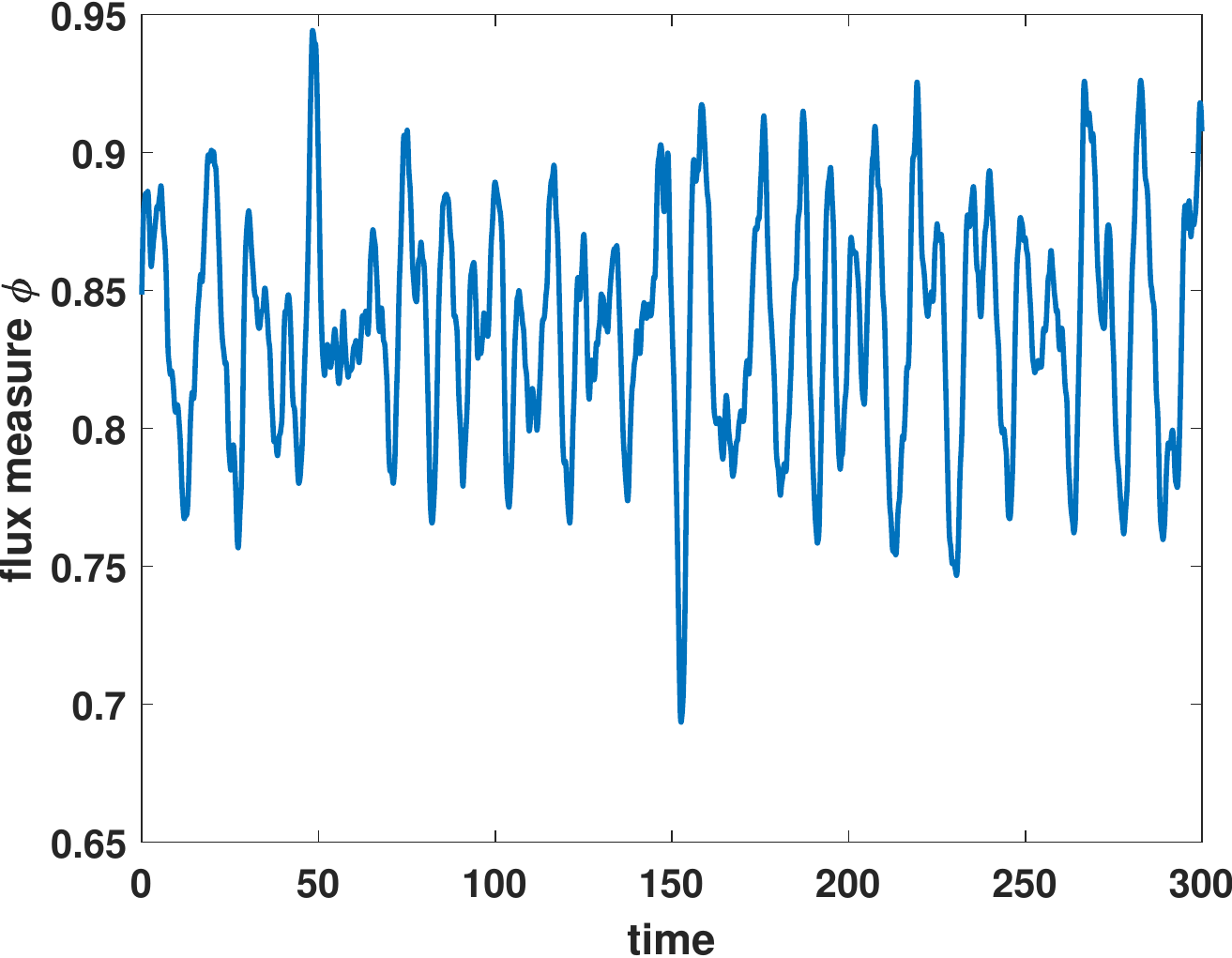}}
\caption{\Cref{fig:fluxm_sweep} shows same data and parameter sweep as \Cref{fig:intro}\textup{(b)} but with different macroscopic measure: output is the value of space-averaged flux $\phi(t)$ at the end of each $300$\,s interval with fixed $\mu$. At this point, we assume that the system has settled at a steady state with fluctuations, as shown in \Cref{fig:fovert}.}
\label{fig:results}
\end{figure}

\Cref{fig:fluxm_sweep} repeats the result shown in
\Cref{fig:intro}(b), but using the instantaneous flux measure
$\phi$ in its $y$-axis. The underlying data is the same in
\cref{fig:intro}(b) such that the same bistability is observed. In
\Cref{fig:fovert}, the time profile of the flux $\phi$ for
$\mu=-1.2$\,m is shown to illustrate the size and time scale of
fluctuations and level of stationarity for a stable steady
flow. Although at this position there is a stable and observable
steady state, the flux exhibits large fluctuations because of the
finite number $N$ of pedestrians. For the value $\mu=-1.2$, for which
\Cref{fig:fovert} is shown, a steady flow of pedestrians moving only
on the right side of the obstacle is observed.  The standard deviation
through the last 20 seconds is 0.05, while the one for the whole time
interval is 0.0418. This measure is used to estimate stationarity.

\subsubsection{Control input through bias force}
\label{subsections:control_input}
For performing continuation on the evacuation scenario with feedback control as proposed in \Cref{sec:1d:SO,sec:1d:combi}, we need to specify how the scalar control input $u$ enters model \eqref{model:xddot}. We also aim to choose an input that is  feasible in the sense that could be applied similarly in a real-life experiment. The motivation comes from traffic control, where signage and traffic light is used to control or direct flows. For our specific scenario, a display showing arrows in front of the obstacle acts as the control input. Pedestrians close to the display see these arrows pointing left or right, where the size of the arrow depends on the control input $u$, while the arrows do not directly influence pedestrians that cannot see the display.

Let us recall the \emph{social force model} \eqref{model:xddot}. We add a linear scalar control input on the particles that induces a behavior similar to a reaction to displaying arrows to pedestrians:
\begin{align}\label{model:xddot+input}
\mathbf{\ddot x}_i&= \mathbf{F}_{\mathrm{al},i} + \sum_{j=1}^N\mathbf{F}_{\mathrm{rep},ij}^\mathrm{ped}
+ \sum_{\mathrm{objects\ }k} \mathbf{F}_{\mathrm{rep},ik}^\mathrm{obj}+
\begin{bmatrix}
b(\mathbf{x}_i)\\ 0
\end{bmatrix}
a u\mbox{,\ where\ }\\
\label{eq:box}
b(\mathbf{x})&=b(x,y)=
\begin{cases}
  1&\mbox{if   $|x-x_{\mathrm{c},b}|\leq x_{\mathrm{wth},b}$, $|y-y_{\mathrm{c},b}| \leq y_{\mathrm{len},b}$}\\
  0&\mbox{otherwise,}
\end{cases}
\end{align}
where $u$ is the scalar time-dependent control input, and $a$ is the scaling of the gains, as introduced for our %general feedback control depending on the extended washout filter 
zero-in-equilibrium feedback control \eqref{combi:ode+control} in \Cref{sec:combi}. 
The parameters $(x_{\mathrm{c},b},y_{\mathrm{c},b})=(\mu,-1.75\,\mathrm{m})$
and $(x_{\mathrm{wth},b},y_{\mathrm{len},b})=(3.3\,\mathrm{m},0.25\,\mathrm{m})$ define a box in
which the factor $b$ in front of the control $u$ has non-zero
support. This box is in front of the obstacle, as illustrated in
\Cref{fig:box} (with red frame). It is the area where the pedestrian
have visual contact with the monitor. The influence of the feedback
control input $u$ is zero outside of this box. The factor $(0,b)$ for
$au$ in \eqref{model:xddot+input} describes that, when a pedestrian is
inside the control input support box, a bias force acts on them
pushing perpendicular to the direction of the corridor (in
$x$-direction, see also \Cref{fig:corridor}). When closing the feedback loop we permit the input $u$
to depend on the output $\phi$, the flux measure defined by
\eqref{flux}.

\subsubsection{Choice of control gains}
\label{sec:gains_choice}
Here we discuss the choice of gains $(a, K_{\mathrm{st},\mu}, K_{\mathrm{st},\phi})$ for the continuation of the particle flow model. Note that the subscript of the last gain (formerly $K_{\mathrm{st},y}$) is now $\phi$ to indicate that the output is $\phi$ as defined in \eqref{flux}. First, we make some assumptions based on experimental evidence, e.g. on the response of the system. Then, the choice of proper gains is possible by combining these assumptions  with properties of the zero-in-equilibrium feedback control%extended washout filter
, see also \cref{sec:1d:combi}.

\paragraph{Input orientation $\sigma$} The blue crosses and the red circles in \cref{fig:intro}(c) show the stable parts of the equilibrium branch
for the flux measure $\phi$ as observed during the parameter sweep, see \eqref{flux} for the definition of $\phi$. Hypothesizing a saddle-node
bifurcations near the transitions and a single unstable connecting
branch between the stable branches, the topology of the figures
implies that the sign of the input orientation $\sigma$ in \eqref{1d:combi:gaincrit3} must be negative at the equilibrium curve at $\mu\approx -1.25$. Thus, $\sigma = -1$.
Recall that $\sigma$ determines
the orientation between equilibrium curve tangent in the
$(\mu,\phi)$-plane and the line $\dot \mu=0$.

\paragraph{Input effect} Next, we determine the sign of the coefficient $w_\mathrm{c}^\tran{}f_u$ appearing in the conditions \eqref{1d:combi:gaincrit3}: By construction, the input $u$ acts  as a bias force, see \eqref{model:xddot+input}, such that
positive $u$ has a direct effect on the flux measure $\phi$ with
positive sign which indicates $w_\mathrm{c}^\tran{}f_u>0$. 

\paragraph{Control gains} With the signs of $\sigma$ and $w_\mathrm{c}^\tran{}f_u$ as determined above, in order to satisfy the conditions in \eqref{1d:combi:gaincrit3} we chose fixed values for the gains $a$ and $K_{\mathrm{st},\phi}$. The parameter $K_{\mathrm{st},\mu}$ is adjusted in each continuation step such that, (a), it is bounded from above, and (b), the line defined in \eqref{1d:combi} is not parallel to the equilibrium curve. Recall that the secant vector $v = (v_{\mu},v_{\phi})$ is the approximation of the tangent to the equilibrium curve. In practice, we chose
\begin{align}
\label{gains_choice}
a&=50,&
K_{\mathrm{st},\phi}&=-0.2,&
|K_{\mathrm{st},\mu}|&= 0.2 \, \frac{v_{\mu}}{v_{\phi}} .
\end{align}
Note that for $K_{\mathrm{st},\mu} = -0.2 v_{\mu}/v_{\phi}$ the line $\dot \mu =0$ is perpendicular to the secant  vector $v$.
However, close to the sign changes of $v_{\phi}$ (extrema in $\phi$ of the equilibrium curve),  \eqref{gains_choice} may lead to a gain $K_{\mathrm{st},\mu}$ with large modulus and possibly positive sign, which can lead to violation of the first condition in \eqref{1d:combi:gaincrit3}.

In these cases, we chose  $K_{\mathrm{st},\mu}=0.2 v_{\mu}/v_{\phi}$ and we confirmed that this choice of gains indeed stabilized the system so that we successfully tracked the equilibrium curve. In practice, we switched the sign of $K_{\mathrm{st},\mu}$ whenever
%we observed whether
during continuation the system was 
driven away from a reference point $(\mu_{\mathrm{ref}}, \phi_\mathrm{ref})$ by more than distance $r_\mathrm{\sigma}=0.2$. For $K_{\mathrm{st},\mu}=0.2 v_{\mu}/v_{\phi}$ the line $\dot \mu = 0$ is not perpendicular to the equilibrium curve. However, the second condition on the gains in \eqref{1d:combi:gaincrit3} is still satisfied.

\section{Control-based continuation of the pedestrian flux}
\label{sec:control-based-ped}
\label{sec:cbped}
In this section, we will test two of the proposed inherently non-invasive feedback control laws, namely the washout filter \eqref{gen:wo:ode} and the  zero-in-equilibrium feedback control %extended washout filter 
introduced in \eqref{combi:ode+control}.
As already briefly discussed in \Cref{sec:1d:combi}, while feedback control through the parameter (described in \Cref{sec:SP,sec:1d:SP}) should be also feasible and stabilizing, we observed that the control would have required large gains $(K_{\mathrm{st},\mu},K_{\mathrm{st},y})$  that violated \Cref{ass:1d} about the presence of only a single slow direction.

Hence, our setup has a control input $u$, given in \eqref{model:xddot+input}, \eqref{eq:box} in \Cref{subsections:control_input}, separate from the also varying bifurcation parameter $\mu$. This is similar to the driven-pendulum experiments performed  by Bureau et al. \cite{bureau2013experimental}, where control was also input through an external force, a real-time controllable magnet.

According to \eqref{1d:SO} in \Cref{sec:1d:SO} for the washout filter we may set the factor $a=1$ in \eqref{model:xddot+input} without loss of generality, and set
\begin{align}
  \label{ytarget}
  u(t)&=K_\mathrm{st}\phi(t)+K_\mathrm{wo}y_\mathrm{wo}(t)\mbox{,}&
  \dot y_\mathrm{wo}&=u(t)\mbox{.}
\end{align}

For zero-in-equilibrium feedback control, we keep $a$ in
\eqref{model:xddot+input} as an adjustable gain (fixed as given in \eqref{gains_choice}) and set, according to
\eqref{1d:combi},
\begin{align}
  \label{ytarget2}
  u(t)&=K_{\mathrm{st},\phi}[\phi(t)-\phi_\mathrm{ref}]+K_{\mathrm{st},\mu}[\mu-\mu_\mathrm{ref}]
  % K_{\mathrm{st}}b(A(f-\tilde{f})+B(\mu-\tilde{\mu}))
  \mbox{,}&
  \dot{\mu}&=u(t)\mbox{.}
\end{align}
During continuation the reference values $(\phi_\mathrm{ref},\mu_\mathrm{ref})$ in \eqref{ytarget2} are defined via a secant prediction. 

\paragraph{Protocol for simulation with feedback control and determination of steady states}
Along the continuation of a branch the simulation runs as a continuous
computational experiment, numerically integrating \eqref{model:xddot+input} with
\eqref{ytarget} or \eqref{ytarget2}. At certain times $t_\mathrm{set}$ some parameters and possibly the gains are set to values determined by the continuation. For the washout-filtered control law \eqref{ytarget} these are the washout filter state $y_\mathrm{wo}(t_\mathrm{set})$ and obstacle $\mu$, while the gains are kept fixed: $(K_\mathrm{st},K_\mathrm{wo})=(-5,\pm0.1)$ with the sign of $K_\mathrm{wo}$ depending on the stability of the branch. For the general zero-in-equilibrium control law \eqref{ytarget2}, the varying parameters and gains are  $(\phi_\mathrm{ref},\mu_\mathrm{ref})$ and $K_{\mathrm{st},\mu}$ while $K_{\mathrm{st},\phi}$ and $a$ are fixed according to \eqref{gains_choice}. After setting these parameters, the simulation is continued until it becomes stationary. At this point %we assume that
a steady state has been reached such that we record it, and determine new parameters $(\phi_\mathrm{ref},\mu_\mathrm{ref})$ and $K_{\mathrm{st},\mu}$.

Our condition determining that a steady state is reached is as follows. When
integrating numerically the system \eqref{model:xddot+input},\,\eqref{ytarget} or \eqref{ytarget2}, we obtain
a sequence ${\phi_0,\phi_1,\phi_2,\ldots,\phi_n}$ of outputs after $n$ integration
steps. We assume that the output measure \emph{$\phi$ is stationary} at time $n$ when the standard
deviation $\std_n$ of last $n_\mathrm{min}$ outputs, $(\phi_{n-n_\mathrm{min}+1},\ldots,\phi_n)$, does not exceed a
given tolerance $\tol_{\std}$. We choose $n_{\min}=200$ which corresponds to $20$
seconds (time step is $0.1$\,s) and tolerance $\tol_{\std}=0.05$. This choice takes into account that we expect from the law of large numbers that the fluctuations will have standard deviation $\sim1/\sqrt{N}$ and is consistent with the observed standard deviation of the time series of the flux
for the stable steady states, as shown in \Cref{fig:fovert}. 

\subsection{Results for control with washout filter}
\label{sec:results:SO}
\begin{figure}[t]
\label{fig:3pieces}
\centering
\includegraphics[scale=0.6]{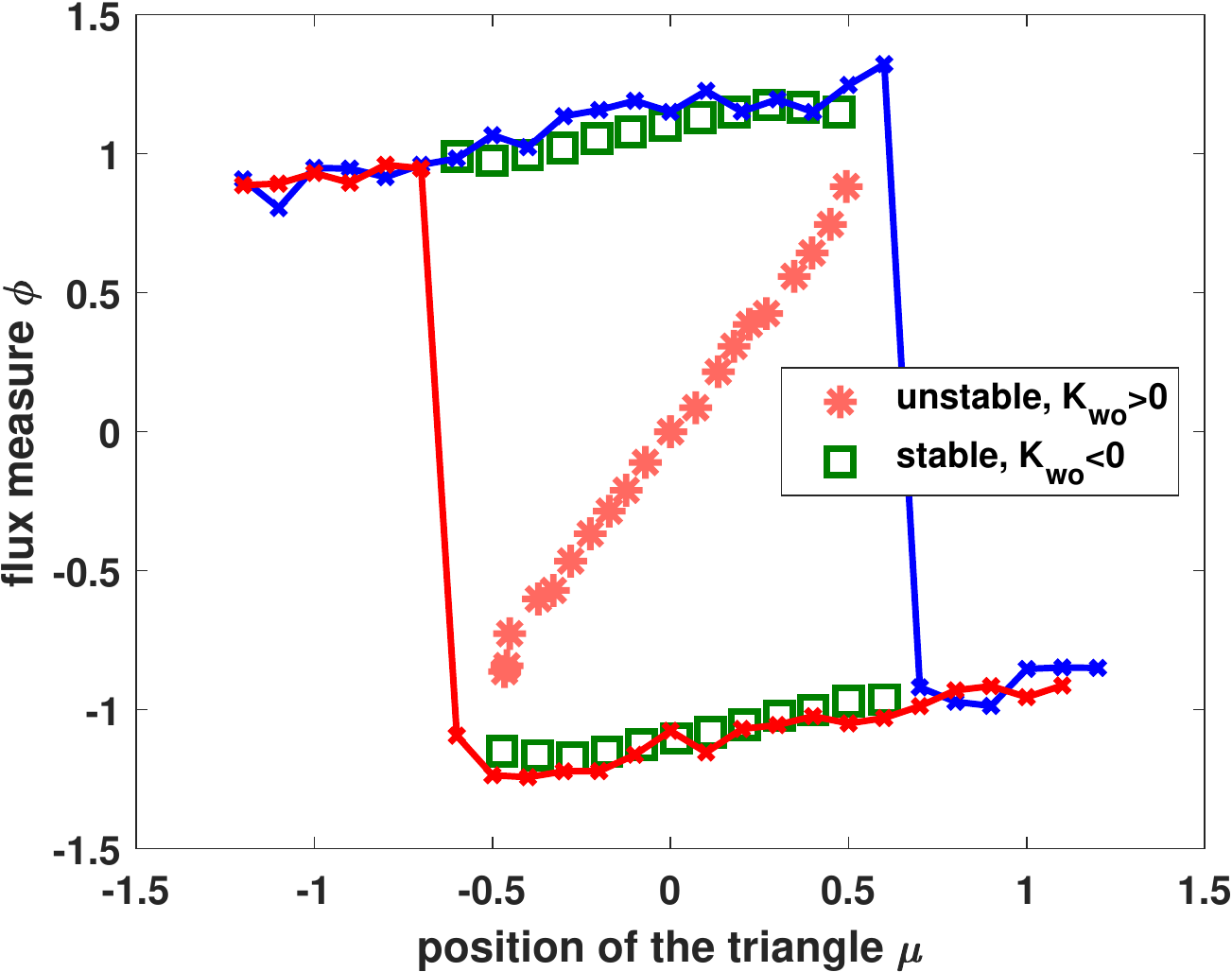}
\caption{Bifurcation diagram of the evacuation scenario obtained by control through washout filters as given in \eqref{ytarget}. The red stars correspond to unstable steady states \textup{(}tracked for  $K_\mathrm{wo}>0$\textup{)} and the green ones correspond to stable steady states ($K_\mathrm{wo}<0$). The blue and red lines  show the observed steady states from sweeping the parameter $\mu$ \textup{(}identical to \Cref{fig:fluxm_sweep}\textup{)}. Control gains: $K_\mathrm{st}=-5$, $K_\mathrm{wo}=\pm0.1$}
\end{figure}
As expected by \Cref{thm:washout} the control using washout filters, \eqref{ytarget}, will be able to track a branch of unstable steady states if one exists. As \Cref{fig:3pieces} demonstrates the control indeed settles to a steady state with fluctuations below tolerance for a range of fixed $\mu$ in the region of unstable states (see the red stars in  \Cref{fig:3pieces}). In this region, as required by our analysis in \Cref{sec:1d:SO}, the gain $K_\mathrm{wo}$ had to be positive. We also recovered the two branches of steady states of the system when using control \eqref{ytarget} and $K_\mathrm{wo}<0$, thus, the completing bifurcation diagram except for small regions near the transitions. 

During the continuation of each branch at each step $j$ (at time
$t_{\mathrm{set},j}$) we used a secant prediction 
$(\mu_{\mathrm{pred},j},\phi_{\mathrm{pred},j})$ with a step size of $\Delta_\mathrm{cont}=0.1$. The parameter $\mu$
is then kept constant at $\mu_{\mathrm{pred},j}$ throughout the
simulation and the initial value for the variable $y_\mathrm{wo}$
equals
$y_\mathrm{wo}(t_{\mathrm{set},j})=-\phi_{\mathrm{pred},j}
K_\mathrm{st}/K_\mathrm{wo}$, such that $u(t_{\mathrm{set},j})=0$  at
the starting time $t_{\mathrm{set},j}$ of the simulation with new
parameters. When the stationarity condition is satisfied we obtain the
output, flux measure $\phi_j$, determining the new approximate fixed
point $(\mu_j,\phi_j)$ with $\mu_j=\mu_{\mathrm{pred},j}$.
The result displayed at \cref{fig:3pieces} is computed in $3$ pieces, as the control \eqref{ytarget} fails to stabilize near the transition points.

The results shown in  \Cref{fig:3pieces} support the hypothesis that there is a unique
unstable branch separating the stable branches in the region of bistability. The diagram suggests that the transition is a saddle-node (or fold) bifurcation but with the control \eqref{ytarget}, using the standard washout filters we are unable to track the branches near these transition points.  
As \Cref{thm:washout} shows, this type of control is singular close to folds in contrast to \eqref{ytarget2}.

\subsection{Results for %control through the extended washout filter
zero-in-equilibrium feedback control}
%\label{sec:res_extended_wo}
\label{sec:res_zie}
\begin{figure}[t]
\centering
\label{fig:results-R}
\includegraphics[scale=0.6]{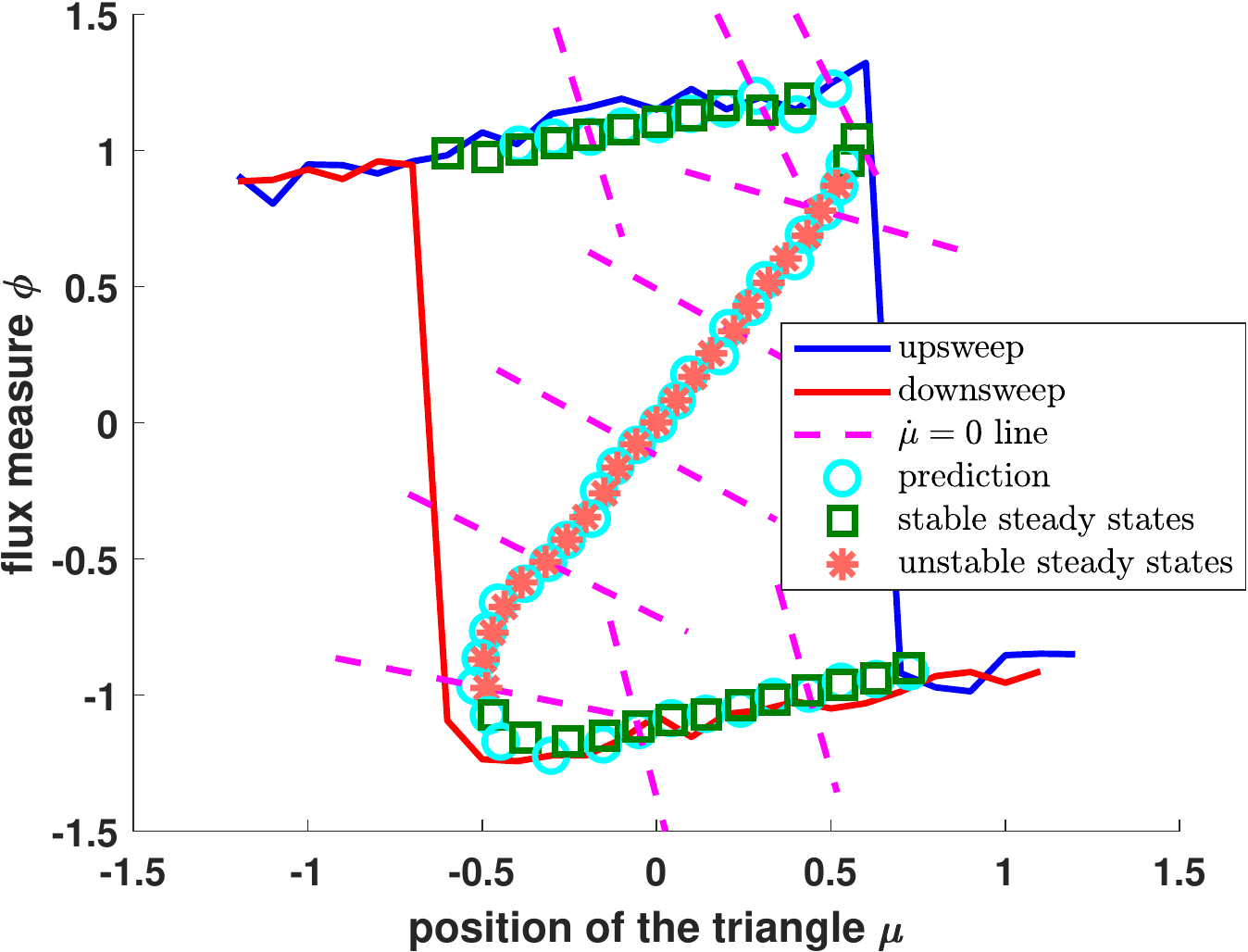}
\caption{Bifurcation diagram of the evacuation scenario obtained by %control through the extended washout filter
zero-in-equilibrium feedback control, \eqref{ytarget2}. A family of equilibria is tracked through two saddle-node bifurcations. A predicted point $(\mu_\mathrm{ref}, \phi_\mathrm{ref})$ is marked with a cyan circle. The fixed gains are $(K_{\mathrm{st},\phi},a) = (-0.2, 50)$. The gain $K_{\mathrm{st},\mu}$ is adjusted for each steady state. Its value is implied by the lines $\dot\mu=0$  \textup{(}only shown for a few selected points in magenta, dashed\textup{)}, see also \eqref{ytarget2}  At each step of the continuation the steady state has to be on its \textup{(}$\dot \mu=0$\textup{)}-line The blue and the red lines show the results of the parameter sweep.
}
\end{figure}
As expected by \Cref{thm:combi}, zero-in-equilibrium feedback control, \eqref{ytarget2}, succeeded in tracking the full bifurcation diagram of the underlying system. The results are demonstrated in \Cref{fig:results-R}.  At step $j$ of the continuation (at time
$t_{\mathrm{set},j}$) the secant prediction for the next steady state is % Therefore, for tracking a new steady state a secant prediction
$(\mu_{\mathrm{pred},j},\phi_{\mathrm{pred},j})$, which enters control law \eqref{ytarget2} as $(\mu_\mathrm{ref}, \phi_\mathrm{ref})$. The step size here is also $\Delta_\mathrm{cont}=0.1$. 
 The slope of the line $\dot\mu=0$ is determined by the gains $K_{\mathrm{st},\phi}$ and $K_{\mathrm{st},\mu}$, and chosen, as described in \Cref{sec:1d:combi} and specified in \eqref{gains_choice},
 %orthogonal to the secant $(v_{\mathrm{sec},\mu},v_{\mathrm{sec},\phi})$
 intersecting the equilibrium curve, and
 with $K_{\mathrm{st},\phi}=-0.2$. The large but fixed choice of gain scaling $a=50$ in \eqref{gains_choice},  with which $u$ enters in \eqref{model:xddot+input}, implies that feedback through the parameter is small in amplitude compared to the feedback through the input shown in \Cref{fig:box}.

\Cref{fig:results-R} also shows the lines $\dot\mu=0$ (dashed magenta lines) for some steps of the continuation to illustrate the adjustment of control gain $K_{\mathrm{st},\mu}$% For illustrative reasons, the target line  $A(f-\tilde{f})+B (\mu-\tilde{\mu})=0$ is plotted (
. The cyan circles are the predictions $(\mu_{\mathrm{pred},j}, \phi_{\mathrm{pred},j})$ for each step $j$. The green squares and red stars are the accepted values for an equilibrium as fluctuations dropped below tolerance.  The green and red colors correspond to stable and unstable equilibria, respectively, as determined by the topology of the branch. For comparison, the blue and the red lines show the results of the parameter sweep (identical to \cref{fig:fovert}). 
  
\begin{figure}[t]
\centering
\label{fig:upper_fold}
\includegraphics[scale=0.6]{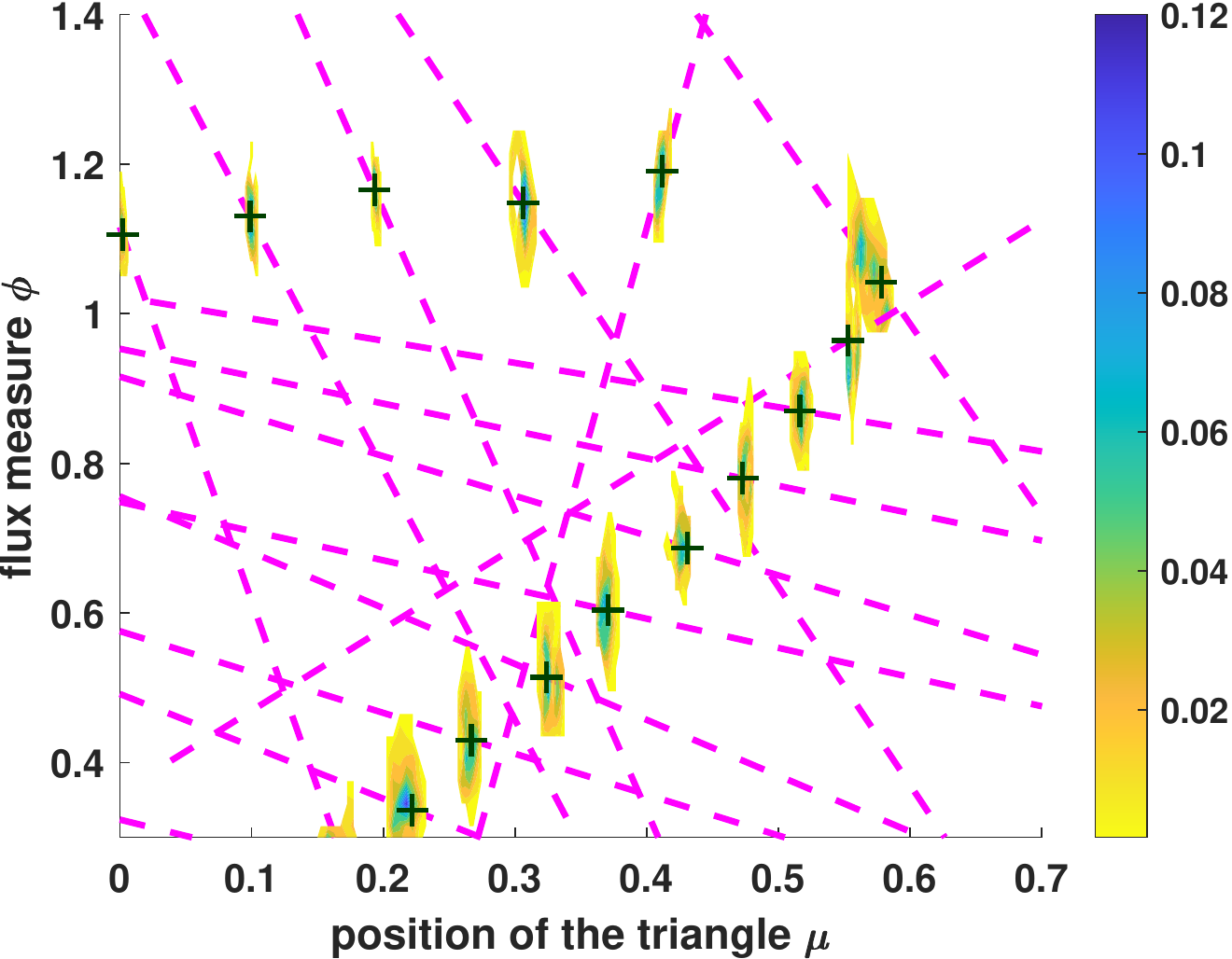}
\caption{Results of %control through extended washout filter
zero-in-equilibrium feedback control close to the fold: evolution of the system in the $(\mu,\phi)$-plane throughout the procedure of tracking steady states close to the upper fold. The last $60$ seconds of the trajectory $(\mu(t),\phi(t))$ are displayed. The accepted fixed points are shown as black crosses.  The magenta dashed lines are the $\dot \mu=0$ lines according to  \eqref{ytarget2}. The color code indicates the percentage of time spend at a position $(\mu,\phi)$.}
\label{fig:colorfold}
\end{figure}

 \Cref{fig:colorfold} shows additional details of the behavior of the system with control law \eqref{ytarget2} close to the upper fold point of the branch. The parameter $\mu$ changes dynamically throughout the simulation. For each fixed point, the dynamical system moves in the $(\mu,\phi)$-plane, eventually approaching a new steady state $(\phi_j,\mu_j)$. The feedback introducing dynamics for $\mu$ permits us to track solution branches around the fold bifurcation, as had been observed by \cite{siettos2004coarse}. \Cref{fig:colorfold} shows the percentage of time spent at every position at the $(\mu,\phi)$-plane during the final $60$ seconds before acceptance of the steady state on a color scale. Away from the fold the parameter $\mu$ is almost constant during transients, apart from the initial adjustment to $\mu_{\mathrm{pred},j}$ (visible in the form of upright horizontally narrow color ellipses in \Cref{fig:colorfold}). However, when approaching the fold point, the variation in the $\mu$ direction becomes larger moving the parameter $\mu$, thus highlighting that control through the parameter plays a role here, even though it is smaller by the factor $a= 50$ compared to the input $u$ in \eqref{model:xddot+input}. The lines $\dot\mu=0$, required to be intersecting the branch for every step  (dashed line in magenta), show where the new steady state should lie along the branch. We  observe that during transients the evolution of the system in the $(\mu,\phi)$-plane does not stay on this line. The stabilization by the feedback control only implies that the system should converge to it. %
% For illustration reasons, only the last 60 seconds of the evolution are displayed.
 %\cref{fig:upper_fold} presents the evolution of the system throughout the calculation of a new steady state while in \cref{fig:upper_fold20}, only the last 20 seconds are plotted.
 \paragraph{\js{Failure of} control through the bifurcation parameter}
\js{After presenting the results of successfully implementing the classical washout filter and the %extended washout filter
zero-in-equilibrium feedback control that was introduced in this paper, we now briefly discuss results of experimentations with applying feedback control only through the parameter (described in \Cref{sec:SP,sec:1d:SP}) %as implemented in  \cite{BartonSieber2013}
to stabilize the system.} 

\begin{figure}[t]
\centering
%\label{fig:upper_fold}
\includegraphics[scale=0.6]{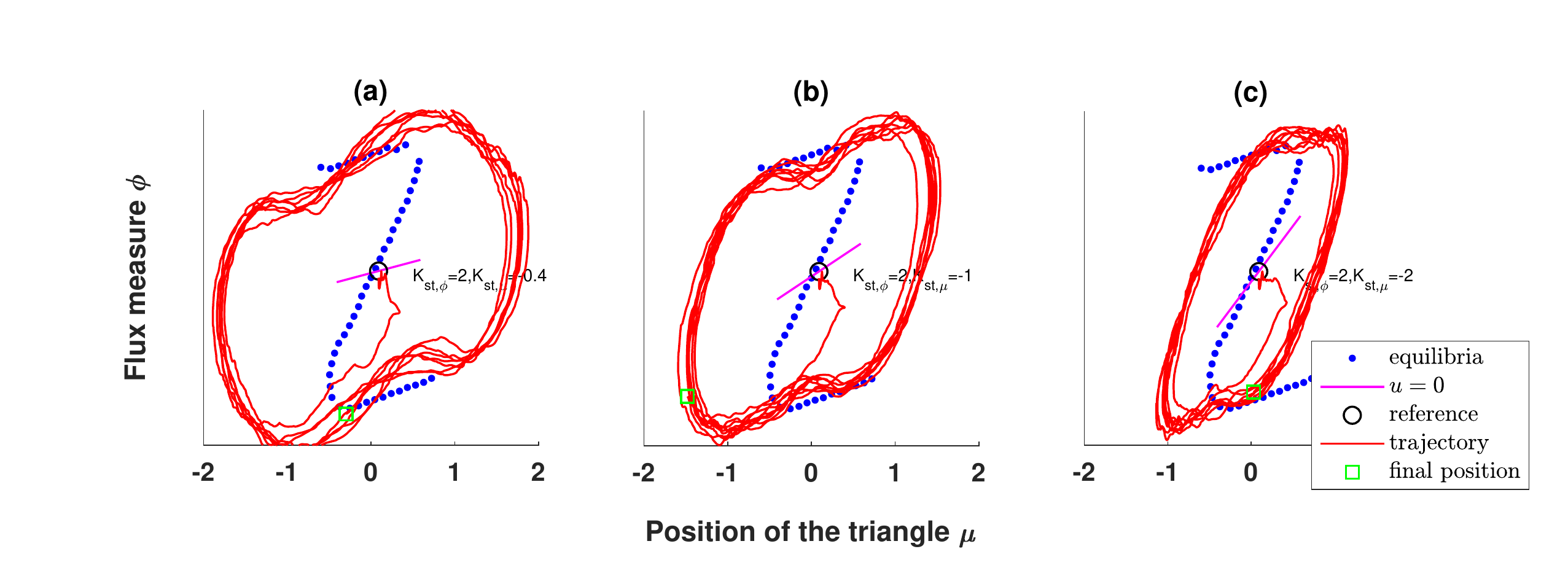}
\caption{\js{Results of %control through extended washout filter
the control through the bifurcation parameter close to the symmetric case of the corridor $(\mu_{\mathrm{pred},j}, \phi_{\mathrm{pred},j})=(0.0895,0.06)$, for different slopes of the line $u=0$. The equilbria as computed by the zero-in-equilbrium feedback control are plotted in blue. The line $u=0$ is plotted in magenda and the reference point with a black circle. The evolution of the system in phase space is plotted in red. For $K_{\mathrm{st},\mu} < 0$ and fixed $K_{\mathrm{st},\phi}=2$ such that {u=0} intersects the equilibrium line with the proper angle (see the discussion in \ref{sec:1d:SO}) , the system oscillates around the control line $u=0$, jumping between the stable branches. It fails to approach the intersection of the line $u=0$ and the equilibrium curve to converge to it as the small fluctuation in the pedestrian flow drives the system away from the unstable equilibrium.}}
\label{fig:ped_line}
\end{figure}

\js{\Cref{fig:ped_line} shows projections of the controlled flow in the $(\mu,\phi)$ plane when control is applied only through the bifurcation parameter $\mu$ near the known unstable equilibrium point $(\mu_\mathrm{eq}(s_0),\phi_\mathrm{eq}(s_0))=(0,0)$ and for various control gains $(K_{\mathrm{st},\phi},K_{\mathrm{st},\mu})$. The topology of the equilibrium branch and our choice $\sigma=-1$ imply that $\partial_s\mu_\mathrm{eq}(s_0)<0$ and $\partial_s\phi_\mathrm{eq}(s_0)<0$. Thus, our criterion for stabilizing gains (assuming that $\phi$ is really governed by a scalar ODE when applying control), \eqref{1d:SP:gaincrit1:final}, requires that (ignoring the $O(\epsilon)$ terms) 
  \begin{align}\label{eq:SPgains_for_ped}
    K_{\mathrm{st},\mu}&<-\lambda_c<0\mbox{,}&
    K_{\mathrm{st},\phi}&> -\frac{\partial_s\mu_\mathrm{eq}(s_0)}{\partial_s\phi_\mathrm{eq}(s_0)}K_{\mathrm{st},\mu}>\frac{\partial_s\mu_\mathrm{eq}(s_0)}{\partial_s\phi_\mathrm{eq}(s_0)}\,\lambda_c>0.
  \end{align}
  These estimates provide a lower bound on the gain $K_{\mathrm{st},\phi}$ depending on the slope of the equilibrium branch and the degree of instability. So, we have to choose $K_{\mathrm{st},\phi}$ positive and sufficiently large, and $K_{\mathrm{st},\mu}$ negative and sufficiently large in modulus (larger than the instability $\lambda_c$). \Cref{fig:ped_line} shows the evolution of the flow when feedback control through parameter is applied with a range of gains satisfying these (only necessary) criteria. % We keep
% This type of control is not suitable for the pedestrian scenario that we studied as the flux measure $\phi$ is quite noisy and fluctuates even on the cases where the pedestrian flow seems rather stable, see for example \cref{fig:intro}. This method did not succeed to track the unstable equilibria of the pedestrian scenarion for all the pairs of  $K_{\mathrm{st},\phi}$ and $K_{\mathrm{st},\mu}$ that we tried. In \cref{fig:ped_line}, we present the evolution of the system in the $(\mu,\phi)$ phase space while applying control through the bifurcation parameter for different slopes of the line $u=0$. Here, w
We fix $K_{\mathrm{st},\phi}=2$ and vary $K_{\mathrm{st},\mu}<0$. %\jsq{we only need gains such that $\{u=0\}$ intersect the equilibrium line the right way and $K_{\mathrm{st},\mu}<0$.} 
The slope of the line $\{u=0\}$ (dashed magenta) indicates the ratio between the gains.  The reference point is $(\mu_{\mathrm{pred},j}, \phi_{\mathrm{pred},j})=(0.0895,0.06)\approx(0,0)$ and is depicted with a black circle. The system shows large-amplitude oscillations around the line $u=0$, jumping between the stable branches of equilibria. The fluctuations of the flux measure $\phi$ cause large and rapid excursions in parameter $\mu$. For larger gains the projection of the flow follows the line $\{u=0\}$ more closely. However, inhererent fluctuations cause correction of the obstacle position that are larger than the corridor width (not shown in \cref{fig:ped_line}). The criteria \eqref{eq:SPgains_for_ped} are necessary, but they are sufficient only under \cref{ass:1d} that the system to be controlled has only a single slow dimension with all others being strongly stable, made in \cref{sec:1dbranch}. \Cref{fig:ped_line} gives evidence that this assumption is violated for gains satisfying \eqref{eq:SPgains_for_ped}. % the slow-fast structure drive the system away from the unstable equilibrium and the control via the bifurcation parameter fails.
}

%We observed that this type of control would have required large gains $(K_{\mathrm{st},\mu},K_{\mathrm{st},\phi})$. These gains violated the assumption $|K_{\mathrm{st},\mu}|\epsilon\ll1$, $|K_{\mathrm{st},\phi}|\epsilon\ll1$ made in the design in \Cref{sec:1d:SP} (see \eqref{1d:SP}) so much that other degrees of freedom of the pedestrian flow were excited,
%adding inertia and stiffness to the system. This violated \Cref{ass:1d} about the presence of only a single slow direction behind our simplifying framework for designing control gains.
%In practice, this control law stabilized the system only for $K_\mathrm{st}=K_{\mathrm{st},\phi}/K_{\mathrm{st},\mu}\ll1$. For these values of gains, the line defined by 
%\begin{align}
%  \mu(t)=\mu_\mathrm{ref}+K_\mathrm{st}\cdot(\phi(t)-\phi_\mathrm{ref})\mbox{,}
%\end{align}
%(see also \eqref{eq:muBS}) results in an almost constant $\mu$ or, equivalently, a line perpendicular to  the $\mu$-axis such that control fails when we track the transition of the system from a stable branch to an unstable branch. 

\section{Discussion and outlook}
\label{sec:conclusions}

In this paper, we give a first outline for a general design principle for inherently non-invasive feedback control laws. These laws can be used to perform bifurcation analysis and track equilibria of dynamical systems that are not given in closed form. Instead we assume that we have an output of an experiment or a macroscopic quantity extracted from a microscopically defined dynamical system. A first result is that we devised a generalization of two well-known non-invasive feedback control laws, namely (a) feedback control with washout filter, which introduces a state observer, and (b) feedback control though the bifurcation parameter, calling our control law zero-in-equilibrium feedback control.

The proposed modification allows for continuous tracking of equilibria in the case of single-parameter studies and, more generally, can be seen as an improvement to the existing methodology as it requires a non-degeneracy condition that fails only at events of codimension $n_x+1$, where $n_x$ is the dimension of the reconstructed state that the feedback control can use. On the contrary, both feedback control with washout filters and feedback control through the bifurcation parameter fail at codimension-one events. Thus, failure becomes generic in single-parameter studies.

We have demonstrated the effectiveness of our feedback control law, the zero-in-equilibrium feedback control, on a particle model for a pedestrian evacuation scenario. For this system, we succeeded in tracking an entire branch of macroscopic equilibria through two macroscopic saddle-node bifurcation points, thus providing evidence that macroscopic bifurcations and unstable states are behind the observed hysteresis and bistability for the particle model. Therefore, our feedback law holds promise also for a successful implementation in a physical experiment, which is much harder to re-initialize after failure of control than computational models.  While we assumed the existence of a low-dimensional model during the design of the feedback control, the results can validate this assumption a-posteriori:  convergence of the controlled system and vanishing control input ensure that the phenomenon observed is natural for the uncontrolled system. This is in contrast to analyzing a derived macroscopic model, as this derivation would be based on assumptions that may be hard to validate in states that cannot be observed in the uncontrolled system due to their dynamical instability or sensitivity.

\paragraph{Outlook for experiments}
As the particle model explicitly included an alignment tendency among
pedestrians moving in similar directions (see \eqref{model:eal}), a
natural next step are controlled physical experiments on pedestrian flows with real humans and
data collection. These experiments will show when such an alignment
tendency and the resulting bistability are really present, which we
expect to depend strongly on the situation. The experiment will
require manufacture of an obstacle that permits real-time adjustment
of its position and an implementation of a bias force for subset of
pedestrians in front of the obstacle.  The analysis of the physical
experiment would in turn lead to 
%a better understanding of
validation of
high-dimensional particle models with multi-scale interactions, and understanding human pedestrian flows in real life.

%\paragraph{Extension to equation-free framework}\jsq{I have included the references under ``equation-free'', so we may not need this anymore.}
%\js{As discussed in \cref{sec:control-based},  equation-free  methods have been developed to perform numerical continuation of a macroscopic property of computational experiments where the dynamics are described at a microscopic scale. An exciting research area is the use of non-invasive feedback control within the equation-free framework so that the control input can 'washout' any model uncertainties. %, in contrast to the typical case where the estimation of jacobians at the macroscopic level is necessary. 
%In particular, we would like to apply the proposed zero-in-equilibrium feedback control within the equation-free framework. 
%Similar work has been done in \cite{siettos2012equation}, where the authors design a washout filter to detect and stabilize the unknown macroscopic steady states of an agent-based model with uncertainty. 
%The authors suggest that control gains can be set through pole placement or discrete linear quadratic regulator techniques. 
%Similarly, the work  in \cite{tsoumanis2012detection}  also performs control-based continuation at a macroscopic level. However, the control gains are set based on geometric arguments without the need to estimate the jacobian of the macroscopic description of the underlying system.
%}

\paragraph{More general inherently non-invasive control
  design principle}
The zero-in-equilibrium feedback control is clearly not the most general possible
formulation for inherently non-invasive feedback control. First,
multiple control inputs may be practical when applying the feedback to
computational experiments, as it may permit simpler design of control
gains. Second, additional integral components (similar to
$x_\mathrm{wo}$), can be used to enforce constraints on the stabilized
equilibrium that either detect or suppress a bifurcation. As a simple
illustrative example, let us consider a limitation of the zero-in-equilibrium feedback control.  The system \eqref{combi:ode+control} is not
controllable in a pitchfork bifurcation point of a system with
reflection symmetry. However, we can use additional integral
components to enforce the symmetry. For example, for a system with, e.g., reflection symmetry $R$ ($Rf(x,\mu,0)=f(Rx,\mu,0)$) we may
formulate a non--invasive feedback law of the form
\begin{align}\label{out:pitchfork}
  \dot x&=f(x,\mu,a_uu+a_\mathrm{wo}x_\mathrm{wo})\mbox{,}&
  \dot \mu&=u\mbox{,}&
  \dot x_\mathrm{wo}&=p^\tran[R-I]x\mbox{}
\end{align}
with $\dim u=\dim x_\mathrm{wo}=1$, weights $a_u,a_\mathrm{wo}\in\R$, and
$p\in\R^{n_x}$ such that the tangent to the asymmetric branch is not
orthogonal to $p$.  One may easily check, that this system is
controllable in the pitchfork bifurcation point $(x,\mu)=(0,0)$ for
the normal form example $f(x,\mu,u)=\mu x+x^3+u$ with symmetry $Rx=-x$
and $p=1$. The key assumption is that the control input is able to
break the symmetry ($p^\tran[R-I]f_u\neq 0$). Thus, feedback control
laws designed for \eqref{out:pitchfork} would permit continuation of
symmetric branches through a pitchfork bifurcation.
This simple illustration shows that more
general bifurcation control is possible. Of particular interest is
control for periodic orbits of autonomous systems. While reduction to
a discrete system through a Poincar{\'e} map should give
straightforward results, this may not be the most practical approach
for a system with large disturbances. Consequently, it would be
interesting to revisit time-delayed feedback control to see if it can
also be formulated in a way similar to the zero-in-equilibrium feedback control.

\section*{Acknowledgments}
Jan Sieber's research was supported by the UK Engineering and Physical
Sciences Research Council (EPSRC) grants EP/N023544/1 and
EP/V04687X/1. Jens Starke's research was supported by the DFG (Deutsche Forschungsgemeinschaft) through the Collaborative Research Center CRC 1270, Grant/Award Number: SFB 1270/1-299150580.

\bibliographystyle{siamplain}
\bibliography{references}
\appendix

\section{Parameter values for model and methods}
\label{app:sec:parameters}
\paragraph{Model parameters}
\Cref{tab:model:parameters} shows all parameters entering the social force  model with alignment, defined in \Crefrange{eq:dir_force}{eq:kappa}, \Cref{tab:flux:parameters}  lists parameters entering pedestrian flux measures $\Phi_\pm$ and $\phi$, defined in \Crefrange{flux}{w}, and \Cref{tab:input:parameters} lists parameters entering the input effect, $b$, given in \eqref{eq:box}.
\begin{table}[ht]
\centering{\small
\begin{tabular}{l l l l} 
 \hline
 Parameter$\phantom{l^{\int}}$ & Symbol & Value & Units \\[0.2ex]
 \hline\hline
 Corridor length$\phantom{l^{\textstyle\int}}$ & $C_\mathrm{len}$ & $20$ & $\mathrm{m}$ \\ 
 %\hline
 Corridor width & $C_\mathrm{wth}$ & $10$ & $\mathrm{m}$ \\
 Number of pedestrians & $N$ & $100$ & - \\
 %\hline
  Triangular obstacle's base side & $L_\mathrm{base}$ & $4$ & $\mathrm{m}$ \\
 %\hline
   Triangular obstacle's leg sides & $L_\mathrm{iso}$ & $3$ & $\mathrm{m}$ \\
 %\hline
   Preferred walking speed & $v_\mathrm{trg}$ & $1.34$ & $\mathrm{ms}^{-1}$ \\
   Target point for walking & $\mathbf{x}_\mathrm{trg}$ & $(0,20)$ & $(\mathrm{m},\mathrm{m})$\\
 %\hline
 Reaction time of pedestrians & $\tau$ & $0.22$ & $\mathrm{s}$ \\[0.2ex]
 %\hline
 Pedestrian-pedestrian repulsion & $V_{\mathrm{rep},ij}^\mathrm{ped}$ & $15$ & $\mathrm{m}^2\mathrm{s}^{-2}$ \\ 
 %\hline
 Pedestrian-pedestrian length scale & $\sigma^\mathrm{ped}$ & $1$ & $\mathrm{m}$ \\
 %\hline
 Pedestrian-Obstacle repulsion & $V_{\mathrm{rep},ij}^\mathrm{obj}$ & $10$ & $\mathrm{m}^2\mathrm{s}^{-2}$ \\
 %\hline
 Pedestrian-Obstacle length scale & $\sigma^\mathrm{obj}$ & $2$ & $\mathrm{m}$ \\[0.2ex]
 %\hline
 Lemming effect parameter & $p_\mathrm{al}$ & $0.75$ & - \\[0.2ex]
 %\hline
Weight function $\kappa$ scaling factor & $\gamma$ & $\exp(1)$ & - \\
 %\hline
Weight function $\kappa$ angle scale & $\beta$ & $0.9$ & - \\
 %\hline
 Weight function $\kappa$ angle scale magnitude & $\alpha$ & $15$ & - \\
 %\hline
Weight function $\kappa$ length scale & $\delta$ & $5$ & $\mathrm{m}$ \\[1ex] 
\hline
$\phantom{l^{\int}}$
\end{tabular}}
\caption{\noindent Model parameters for \Crefrange{eq:dir_force}{eq:kappa}, see
  also \Cref{fig:corridor} for meaning of geometry parameters for
  corridor and obstacle.}
\label{tab:model:parameters}
\end{table}
\begin{table}[ht]
\centering{\small
\begin{tabular}{l l l l} 
 \hline
 Parameter$\phantom{l^{\int}}$ & Symbol & Value & Units \\[0.2ex]
 \hline\hline
 %\hline
 ($\Phi_{\pm}$) half-length (in $y$) of box in \Cref{fig:pedweight}$\phantom{l^{\textstyle\int}}$ & $y_{\mathrm{len},\Phi}$ & $0.5$ & $\mathrm{m}$ \\
 ($\Phi_{\pm}$) center (in $y$) of box in \Cref{fig:pedweight} & $y_{\mathrm{c},\Phi}$ & $2.75^{[*]}$ & $\mathrm{m}$ \\
  ($\Phi_{\pm}$) length of interval for time averaging  $\phantom{l^{\textstyle\int}}$ & $\tau_\mathrm{max}$ & $10$ & $\mathrm{sec}$ \\
 % ($\Phi_{\pm}$) separating line (in x dir) of the box for counting pedestrians  %$\phantom{l^{\textstyle\int}}$ & $y_\mathrm{cnt,len}$ & $0.5$ & $\mathrm{m}$ \\
($\phi$) length scale of $w_\mathrm{pos}$ $\phantom{l^{\textstyle\int}}$ & $d$ & $4$ & $\mathrm{m}$ \\
 %\hline
($\phi$) scaling factor/dimension of flux & $\eta$ & $1/12$ & $\mathrm{s}/\mathrm{m}^3$ \\
 %\hline
($\phi$)  distance of extrema for $w_\mathrm{pos}$ from middle of corridor & $x_{\mathrm{c},\phi}$ & $C_\mathrm{wth}/2=5$ & $\mathrm{m}$ \\
 %\hline
($\phi$)  location of extrema for $w_\mathrm{pos}$ along corridor & $y_{\mathrm{c},\phi}$ & $0$ & $\mathrm{m}$ \\[0.5ex]
%\jsq{$\Phi$ parameters still missing} \\[1ex] 
\hline
$\phantom{l^{\int}}$
\end{tabular}}
\caption{\noindent Parameters, defining output measures, time-averaged fluxes
  $\Phi_\pm$ and space-averages flux $\phi$
  in \eqref{flux},\,\eqref{w}, see also \Cref{fig:pedweight} for
  illustration of graph for $w_\mathrm{pos}$ and box location for
  $\Phi_\pm$. $[*]$ The value of $y_{c,\Phi}$ is determined as $\sqrt{\smash[b]{L_\mathrm{iso}^2-L_\mathrm{base}^2/4}}+0.5 = \sqrt 5+0.5\approx 2.75$}
\label{tab:flux:parameters}
\end{table}
\begin{table}[ht]
\centering{\small
\begin{tabular}{l l l l} 
 \hline
 Parameter$\phantom{l^{\int}}$ & Symbol & Value & Units \\[0.2ex]
 \hline\hline
 %\hline
half-length (in $y$) of  box for feedback input $\phantom{l^{\textstyle\int}}$& $y_{\mathrm{len},b}$ & $0.25$ & $\mathrm{m}$ \\
half-width (in $x$) of box for feedback input & $x_{\mathrm{wth},b}$ & $C_\mathrm{wth}/3=3.3$ & $\mathrm{m}$ \\
center (in $y$) of  box for feedback input & $y_{\mathrm{c},b}$ & $-1.75$ & $\mathrm{m}$ \\
center (in $x$) of  box for feedback input & $x_{\mathrm{c},b}$ & $\mu$ & $\mathrm{m}$ \\
\hline
$\phantom{l^{\int}}$
\end{tabular}}
\caption{\noindent Parameters, defining input support box for $b$,
  where the feedback control $u$ can act, defined in
  \eqref{model:xddot+input},\,\eqref{eq:box}, see also \Cref{fig:box}
  for illustration of the location of the box relative to corridor and
  obstacle. } 
\label{tab:input:parameters}
\end{table}

\paragraph{Numerical integration details for simulation of social force model with alignment}
%\label{app:sec:numerical_integration}
For all simulations, we use the \textit{MATLAB} ode45-solver with a fixed step size of 0.1 sec. Throughout the simulation, $N = 100$ pedestrians/particles were always inside the corridor. To this end, an event function was implemented. This functions detects when a particle has left the corridor or, equivalently, when for a particle $i$ with  positions $(x_i,y_i)$ it holds that $y_i>C_\mathrm{len}/2$. In that event, the differential equations that correspond to  particle $i$ are removed from the integrated system and a new particle $i$ is injected at the other side of the corridor. The new particle enters with $y_i=-C_\mathrm{len}/2$ and with a vertical position $x_i$ uniformly random distributed around the center line of the corridor within the interval $[-0.5,0.5]$\,m. The initial velocity is $(v_\mathrm{trg},0)$. 
\end{document}